\documentclass[10pt,letterpaper]{article}

\usepackage[utf8]{inputenc}
\usepackage{url}
\usepackage{epsfig}
\usepackage{graphicx}
\usepackage{subfigure}
\usepackage{mathtools}
\usepackage{algpseudocode}
\usepackage{algorithm}
\usepackage{amsmath}
\usepackage{relsize}
\usepackage{mathrsfs}
\usepackage{bm}
\usepackage{amsfonts}
\usepackage{color}
\usepackage{xcolor}
\usepackage{verbatim}
\usepackage{subfig}

\usepackage[right]{lineno}

\newcommand{\bk}{\mathbf{k}}
\newcommand{\bx}{\mathbf{x}}
\newcommand{\bX}{\mathbf{X}}
\newcommand{\bff}{\mathbf{f}}

\newcommand{\bY}{\mathbf{Y}}
\newcommand{\bI}{\mathbf{I}}
\newcommand{\bK}{\mathbf{K}}

\newcommand{\btheta}{\boldsymbol{\theta}}
\newcommand{\bpsi}{\boldsymbol{\psi}}

\newcommand{\Q}{\boldsymbol{Q}}

\newcommand{\E}{\mathbb{E}}
\newcommand{\R}{\mathbb{R}}
\newcommand{\GP}{\operatorname{GP}}

\newcommand{\calN}{\mathcal{N}}

\newcommand{\sref}[1]{Sec.~\ref{sec:#1}}
\newcommand{\qref}[1]{Eq.~(\ref{eqn:#1})}

\newcommand{\fref}[1]{Fig.~\ref{fig:#1}}

\DeclareMathOperator\erf{erf}

\definecolor{Gray}{gray}{.25}

\begin{document}
\vspace*{0.1in}

\begin{center}
{\Large
\textbf\newline{Bayesian Optimal Design of Experiments For Inferring 
The Statistical Expectation Of A Black-Box Function}
}
\newline
\\
Piyush Pandita\textsuperscript{*},
Ilias Bilionis\textsuperscript{},
Jitesh Panchal\textsuperscript{}
\\
\bigskip
School of Mechanical Engineering, Purdue University, West Lafayette, Indiana 47907
\\
\bigskip
* ppandit@purdue.edu

\end{center}

\begin{abstract}
\label{sec:abs}
{\it 
Bayesian optimal design of experiments (BODE) has been successful in acquiring information about a quantity of interest (QoI) which depends on a black-box function.
BODE is characterized by sequentially querying the function at specific designs selected by an infill-sampling criterion.
However, most current BODE methods operate in specific contexts like optimization, or learning a universal representation of the black-box function.
The objective of this paper is to design a BODE for estimating the statistical expectation of a physical response surface.
This QoI is omnipresent in uncertainty propagation and design under uncertainty problems.
Our hypothesis is that an optimal BODE should be maximizing the expected information gain in the QoI. 
We represent the information gain from a hypothetical experiment as the Kullback-Liebler (KL) divergence between the prior and the posterior probability distributions of the QoI.
The prior distribution of the QoI is conditioned on the observed data and the posterior distribution of the QoI is conditioned on the observed data and a hypothetical experiment. 
The main contribution of this paper is the derivation of a semi-analytic mathematical formula for the expected information gain about the statistical expectation of a physical response.
The developed BODE is validated on synthetic functions with varying number of input-dimensions. 
We demonstrate the performance of the methodology on a steel wire manufacturing problem. 
}
\end{abstract}

\section{Introduction}
\label{sec:intro}
Engineering problems require either computationally intensive computer codes~\cite{sacks1989} or expensive physical experiments~\cite{flournoy1993}. 
With insufficient information about the analytic dependence of the physical response on the design parameters or experimental conditions, the engineer needs scores of physical response evaluations to make decisions with confidence.
To overcome this issue, researchers have developed design of experiments (DOE) techniques that attempt to select the maximally informative physical response evaluations within a given budget~\cite{eriksson2000,anderson2000,alexanderian2014optimal}.
Classical DOE techniques generate a single batch design~\cite{montgomery2017design} and, thus, they face several shortcoming in case of functions with local features, e.g., discontinuities, or sharp non-linearities~\cite{chaloner1995bayesian}.
Sometimes the DOE obtained can be equally spaced when the context requires more samples from certain regions of the domain.
Such scenarios require a sequential DOE (SDOE) approach.

SDOE uses past observations to decide the next evaluation point~\cite{chernoff1959sequential,robbins1985some}.
Over the past two decades, SDOE has been used in several applications spanning both physical experiments~\cite{flournoy1993,havinga2013sequential,havinga2017sequential,alrefae2018process,saviers2017scaled} and computer simulations~\cite{schonlau1997computer,simpson2001sampling,gramacy2015sequential,huan2010accelerated}.
One of the most theoretically sound SDOEs is Bayesian optimal design of experiments (BODE).
Under BODE, one models the physical response using a statistical surrogate and selects the next evaluation point by attempting to maximize the expected value of information.
The newly acquired information is used to condition one's belief about the physical response using Bayes' rule. 
The process is repeated until the marginal value of information is negative.
The exact definition of the value of information depends on one's goals.
For example, one could be interested in optimizing an objective~\cite{locatelli1997bayesian,jones1998,gaul2014modified,huang2006global,lizotte2008,frazier2008knowledge,mockus2012bayesian,arendt2013objective,huan2014gradient,lam2016,marco2016automatic,kristensen2017adaptive,christen2011advances}, learning an accurate representation of the physical response~\cite{mackay1992information,krause2008near,stroh2017sequential,beck2016sequential,gramacy2009adaptive,terejanu2012bayesian} or estimating the probability of a rare event~\cite{mohamad2017direct,mohamad2018sequential}. 

Instead of the value of information, several BODE approaches attempt to maximize the information gain about a quantity of interest (QoI).
The information gain can be quantified through the Kullback-Leibler divergence (KLD)~\cite{kullback1951information,mckay2000} (also known as relative entropy).
Over the years, KLD has been used to quantify information gain~\cite{tsilifis2017efficient} about the objective function, from a hypothetical experiment (an untried design).
The efficacy of the KLD has been extended and demonstrated on various applications including the sensor placement problem~\cite{nath2017sensor,huan2014gradient}, surrogate modeling~\cite{yan2018gaussian,choi2004polynomial,hadigol2017least}, learning missing parameters~\cite{terejanu2013bayesian}, optimizing an expensive physical response~\cite{hennig2012entropy}, calibrating a physical model~\cite{guestrin2005near,huan2013simulation}, reliability design~\cite{picheny2010adaptive}, efficient design space exploration~\cite{liu2016adaptive}, probabilistic sensitivity analysis~\cite{liu2006relative}.

Despite the significant progress, deriving BODE methods for new objectives remains a non-trivial task.
In particular, there are no BODE methods for efficiently propagating input uncertainties through a physical response surface, e.g., estimating the statistical expectation, the variance, or higher order statistics of a physical quantity of interest.
Uncertainty propagation is particularly important for characterizing the robustness of a simulation/experiment and, thus, being able to do it efficiently is essential for robust design.
To address this need, the \emph{objective} of this paper is to develop a BODE methodology for estimating the statistical expectation of the physical response.
The technical details of our approach are as follows. 
Much like the majority of the work in BODE, we use Gaussian process (GP) surrogates to emulate the physical response~\cite{o2006bayesian}.
The expected information gain from a hypothetical experiment is defined to be the KLD between one's \emph{prior} and \emph{posterior} probability densities on the statistical expectation of the physical response.
To derive analytical expressions of the prior and the posterior of this quantity of interest, we use the standard expressions for the mean and covariance of a GP conditioned on data. 
The EKLD of the statistical expectation of the physical response comes out to be an analytically tractable function which alleviates the need for sample averaging.

In summary, our main contributions are as follows: (a) The derivation of semi-analytical expressions for the expected information gain in one's state of knowledge about the statistical expectation of an expensive-to-evaluate physical response; (b) The numerical investigation of the performance of the resulting BODE using synthetic examples; (c) Numerical comparisons to uncertainty sampling; (d) The application of the new scheme to solve an uncertainty propagation problem involving a steel wire manufacturing process simulated using finite elements; and (e) A freely available Python implementation of our methodology\footnote{\url{https://github.com/piyushpandita92/bode}}.

The rest of the paper is organized as follows: \sref{metho} describes in detail
the methodology used, including GP regression \sref{gpr} and the EKLD \sref{ekld}.
The results obtained for three synthetic test problems have been presented in \sref{results}. 
We compare the above proposed BODE methodology with uncertainty sampling which is a common design of experiments method used in practical engineering scenarios in \sref{comp}.
The steel wire manufacturing problem is briefly explained and treated with the proposed methodology in \sref{wire}.
We summarize the nuances of the methodology including its weaknesses and comment on future research directions in \sref{conc}.

\section{Methodology}
\label{sec:metho}
Throughout the paper we represent the various elements of our state of knowledge and objective as follows:
\begin{enumerate}
\item{$\bX_{n}$ are the $n$ designs at which the simulation/experiment has been conducted, i.e., $\bX_n = \{\bx_{1},\cdots,\bx_n\}$}.
\item{$\bY_{n}$ are the values of the physical response at the corresponding $n$ designs, i.e., $\bY_n = \{y_{1}, \cdots, y_{n}\}$}.
\item{Collectively, we represent all observed data by $\mathbf{D}_{n}$ = \{$\bX_{n}$, $\bY_{n}$\}.}
\item{A hypothetical untried design is denoted by $\tilde{\bx}$.}
\item{A hypothetical observation at $\tilde{\bx}$ is denoted by $\tilde{y}$.}
\end{enumerate}
Let $\bx$ be a random variable with probability density function (PDF) $p(\bx)$.
Without loss of generality, we will assume that $p(\bx)$ is the uniform PDF supported on the $d$-hypercube $\mathcal{X}=\times_{k=1}^{d}[0, 1]$.
The true physical response $f$ is assumed to be a squared integrable function of $\mathbf{x}\in\mathcal{X}$, i.e., $f\in\mathcal{L}^{2}(\mathcal{X})$, where
\begin{equation}
    \label{eqn:def_func}
    \mathcal{L}^{2}(\mathcal{X}) = \left\{f:\mathcal{X} \rightarrow \R \middle|\int_{\mathcal{X}}f^{2}(\bx)p(\bx)d\bx < \infty\right\}.
\end{equation}
The QoI $\Q$ that we want to discover through the sequential design of experiments is the statistical expectation of the physical response.
Mathematically, 
\begin{eqnarray}
    \label{eqn:qoi_ex}
    \Q[f] = \int_{\mathcal{X}} f(\bx)p(\bx)d\bx.
\end{eqnarray}
This QoI is a bounded linear functional, an observation that leads to analytical progress.
At each stage of the SDOE, we will update our beliefs about $\Q$ in a Bayesian way, quantifying the epistemic uncertainty induced by limited data at the same time. 
We will select the new experiment by maximizing the expected information gain for $\Q$.

\subsection{Surrogate modeling}
\label{sec:gpr}
GP regression is a very popular non-parametric Bayesian regression technique.
It allows one to express their prior beliefs about the underlying response surface, but it also quantifies epistemic uncertainty induced by limited observations. 
Here, we describe the GP regression very briefly.
More details can be found in~\cite{rasmussen2006}.

\subsubsection{Prior Gaussian process}
\label{sec:prior_gp}
We model our prior beliefs about the physical response using a zero mean GP. 
The covariance function is defined by a radial basis function (RBF), also known as squared exponential. 
Mathematically,
\begin{equation}
    \label{eq:prior_gp}
    f \sim \GP(0, k),
\end{equation}
where
\begin{equation}
    \label{eqn:cov_kern}
    k(\bx,\bx') = k(\bx,\bx';\bpsi) = {s^2}\exp \left\{ { - \frac{1}{2}\sum\limits_{j = 1}^d}{\frac{{{{({x_j} - {x_j}')}^2}}}{{\ell_j^2}}} \right\}.
\end{equation}
The covariance function defined in \qref{cov_kern} encodes our prior beliefs about the smoothness and magnitude of the response. 
The symbol $\ell_{j}>0$ in  \qref{cov_kern} is the lengthscale of the $j$-dimension of the input space.
This parameter quantifies the correlation between the function values at two different inputs.
The $s^2$ in \qref{cov_kern} is the signal strength of the GP.
It incorporates the scale of the response.
These parameters are the hyper-parameters of the covariance function and we will denote them by $\bpsi$, i.e., $\bpsi = \{s^{2}, \ell_{1}, \cdots, \ell_{d}\}$. 
A nonzero mean function can always be included with only minor modifications in what follows.

\subsubsection{The data likelihood}
\label{sec:like_gp}
The likelihood of the data $\bY_n$ a multivariate Gaussian.
The mean vector of this Gaussian distribution is the vector of function output values $\bff_{n}=\left\{f(\bx_1),\dots,f(\bx_n)\right\}$ at observed designs.
The covariance matrix can be computed using the structure defined in \qref{cov_kern}. 
The observations are assumed to be contaminated with Gaussian noise with variance $\sigma^{2}$.
This noise variance is could will be very small relative to the signal strength in the case of computer simulation design.
We augment the vector of hyper-parameters to include this additional parameter to get $\btheta = \{\bpsi, \sigma^{2}\}$.
Mathethematically, the likelihood of the observed data is:
\begin{equation}
    \label{eqn:like_gp}
    p(\bY_{n}|\bX_{n}, \btheta) = \mathcal{N}(\bY_{n} | \bff_{n}, \bK_{n} + \sigma^{2}\bI_{n}) ,
\end{equation}
where $\bK_{n}$ is a $n\times n$ covariance matrix defined according to \qref{cov_kern}, i.e., $K_{nij} = k(\bx_i, \bx_j)$.

\subsubsection{Training the hyper-parameters}
\label{sec:training}
Typically, the hyper-parameter values are fitted to the observed data by maximizing the likelihood of \qref{like_gp}.
However, this process may result in overfitting which is particularly problematic in the context of SDOE.
In this work, we opt for a fully Bayesian treatment~\cite{gelman2014Bayesian} which is more robust.
We assume that the hyperparameters are a priori independent following an exponential prior distribution on the lengthscales and Gamma prior distribution on the signal strength.
Since we do not treat noisy problems in this work, we fix the variance of the likelihood probability to 1e-6 which is a reasonably small value. 
Bayes' rule allows  yields the hyperparameter posterior:
\begin{equation}
    \label{eqn:hyper_post}
    p(\btheta | \mathbf{D}_n) \propto p(\bY_{n}|\bX_{n}, \bpsi) p(\bpsi).
\end{equation}
Here, we employ a \emph{parallel-chain} Markov chain Monte Carlo (MCMC) algorithm with an \emph{affine invariance} sampler to sample from the posterior.
More details on the inner workings of the MCMC algorithm can be found in ~\cite{goodman2010ensemble}.
The code for this MCMC algorithm is available online.\footnote{\url{https://github.com/dfm/emcee}}

\subsubsection{Making predictions}
\label{sec:post_gp}
Conditioned on the hyperparameters, our state of knowledge about $f$ is also characterized by a GP:
\begin{equation}
    \label{eqn:post_gp}
    f|\mathbf{D}_n,\btheta \sim \GP(f|m_n, k_n),
\end{equation}
where
\begin{equation}
    \label{eqn:posterior_mean}
    m_n(\bx) = \left(\bk_n(\bx)\right)^{T}\left(\bK_n + \sigma^2\bI_n\right)^{-1}\bY_{n},
\end{equation}
with  
\begin{equation}
    \bm{\alpha}_n = (\bK_n + \sigma^2\bI_n)^{-1}\bY_{n},
    \label{eqn:alpha_posterior_mean}
\end{equation}
is the \emph{posterior mean} function, and
\begin{equation}
    \label{eqn:predictive_covariance}
    k_n(\bx,\bx') = k(\bx, \bx') - \left(\bk_{n}(\bx)\right)^{T}\left(\bK_n + \sigma^2\bI_n\right)^{-1} \bk_n(\bx'),
\end{equation}
with $\bk_n(\bx) = \left(k(\bx,\bx_1), \dots, k(\bx,\bx_n)\right)^T$, is the \emph{posterior covariance} function.
In particular, at an untried design point $\tilde{\bx}$ the point-predictive posterior probability density  of $\tilde{y} = f(\tilde{\bx})$ conditioned on the hyperparameters is:
\begin{equation}
    \label{eqn:point_predictive}
    p(\tilde{y}|\tilde{\bx}, \mathbf{D}_n, \btheta) = \calN\left(\tilde{y}\middle|m_n(\tilde{\bx};\btheta), \sigma_n^2(\tilde{\bx};\btheta)\right),
\end{equation}
where $\sigma_n^2(\tilde{\bx};\btheta) = k_n(\tilde{\bx},\tilde{\bx};\btheta)$.
Finally, the \emph{point-predictive posterior} PDF of $\tilde{y} = f(\tilde{\bx})$ is:
\begin{equation}
    \label{eqn:posterior_predictive}
    p(\tilde{y} | \tilde{\bx}, \mathbf{D}_n) = \int p(\tilde{y}|\tilde{\bx}, \mathbf{D}_n, \btheta) p(\btheta|\mathbf{D}_n)d\btheta. 
\end{equation}
The latter is, of course, not analytically available, but one can derive sampling average approximations using the MCMC samples from $p(\btheta|\mathbf{D}_n)$.

\subsection{Sequential design of experiments using the expected information gain}
\label{sec:ekld}
Given $\mathbf{D}_n$ observations, our state of knowledge about the QoI $Q[f]$ is given by:
\begin{equation}
    \label{eqn:current_Q}
    p(\Q|\btheta,\mathbf{D}_n) = \mathbb{E}\left[\delta\left(\Q - \Q[f]\right)\middle|\btheta,\mathbf{D}_n \right],
\end{equation}
where the expectation is over the function space measure defined by the posterior GP, see Eq.~(\ref{eqn:post_gp}).
The uncertainty in $p(\Q|\mathbf{D}_n)$ represents our epistemic uncertainty induced by the limited number of observations in $\mathbf{D}_n$.
Now suppose that we did an experiment at $\tilde{\bx}$ and observed the output $\tilde{y}$.
The posterior GP measure would become $p(\Q|\mathbf{D}_n, \tilde{\bx}, \tilde{y})$ and, thus, our state of knowledge about $\Q$ would be:
\begin{equation}
    \label{eqn:hypothetical_Q}
    p(\Q|\btheta,\mathbf{D}_n, \tilde{\bx}, \tilde{y}) = \mathbb{E}\left[\delta\left(\Q - \Q[f]\right)\middle | \btheta,\mathbf{D}_n, \tilde{\bx}, \tilde{y}\right].
\end{equation}
According to information theory, the information gained through the hypothetical experiment $(\tilde{\bx}, \tilde{y})$ conditioned on the hyperparameters, say $G(\tilde{\bx}, \tilde{y};\btheta)$ is given by the KLD between $p(\Q|\btheta,\mathbf{D}_n, \tilde{\bx}, \tilde{y}))$ and $p(\Q|\btheta,\mathbf{D}_n)$.
Mathematically, it is:
\begin{equation}
    \label{eqn:kld}
    G(\tilde{\bx}, \tilde{y};\btheta) = \int_{-\infty}^{\infty}p(\Q|\btheta,\mathbf{D}_n, \tilde{\bx}, \tilde{y}))\log\frac{p(\Q|\btheta,\mathbf{D}_n, \tilde{\bx}, \tilde{y}))}{p(\Q|\btheta,\mathbf{D}_n)}d\Q.
\end{equation}
The expected information gain of the hypothetical experiment, say $G(\tilde{\bx})$, is obtained by taking the expectation of $G(\tilde{\bx},\tilde{y})$ over our current state of knowledge.
Specifically,
\begin{equation}
    \label{eqn:expected_info_gain}
    G(\tilde{\bx}) = \int G(\tilde{\bx}, \tilde{y};\btheta)p(\tilde{y}|\btheta,\tilde{\bx}, \mathbf{D}_{n})p(\btheta|\mathbf{D}_n)d\tilde{y}d\btheta.
\end{equation}
We pick the next experiment by solving:
\begin{equation}
    \label{eqn:iaf_def}
    \bx_{n+1} = \arg\max_{\tilde{\bx}} G(\tilde{\bx}).
\end{equation}
In the rest of this section, we derive analytical approximations of $p(\Q|\btheta,\mathbf{D}_n)$ (\sref{current_state_of_knowledge}), $p(\Q|\btheta,\mathbf{D}_n, \tilde{\bx}, \tilde{y}))$ (\sref{hypothetical_Q}), $G(\tilde{\bx},\tilde{y};\btheta)$ (\sref{information_gain}), and a sampling average approximation for $G(\tilde{\bx})$ (\sref{information_gain}).

\subsubsection{Quantification of the current state of knowledge about QoI}
\label{sec:current_state_of_knowledge}
We now derive an analytical approximation of our current state of knowledge about the QoI, i.e., $p(\Q|\btheta,\mathbf{D}_n)$.
Since the QoI $\Q$, \qref{qoi_ex}, is linear and the point predictive PDF of $y=f(x)$ is Gaussian, \qref{point_predictive}, $p(\Q|\btheta,\mathbf{D}_n)$ is Gaussian.
In particular, it is easy to show that:
\begin{equation}
    \label{eqn:dist_1}
    p(Q|\btheta,\mathbf{D}_n)=\mathcal{N}\left(Q\middle|\mu_1, \sigma_1^2\right).
\end{equation}
The mean $\mu_1$ is given by:
\begin{equation}
    \label{eqn:mu_q}
    \begin{array}{ccc}
    \mu_{1}&:=& \E[\Q|\btheta, D_{n}] \\
    &=& \E\left[\int_{\mathcal{X}} f(\bx) p(x)d\bx \middle |\btheta, D_{n}\right]\\
    &=&\int_{\mathcal{X}}\E\left[f(\bx)\middle |\btheta, D_{n}\right]p(\bx)d\bx \\
    &=&\int_{\mathcal{X}} m_n(\bx)p(\bx)d\bx  \\
    &=& \bm{\epsilon}_n^T\bm{\alpha}_n, 
    \end{array}
\end{equation}
where $\bm{\alpha}_n$ is defined in \qref{alpha_posterior_mean} and each component of $\bm{\epsilon}_n\in\mathbb{R}^n$ is given by:
\begin{equation}
    \label{eqn:epsilon_n}
    \begin{array}{ccc}
    \epsilon_{ni} &=& \epsilon(\bx_i)\\
    &:=& \int_{\mathcal{X}} k(\bx_i, \bx) p(\bx)d\bx \\ &=& s^2\left(\frac{\pi}{2}\right)^{\frac{d}{2}}\prod_{k=1}^d\left\{\ell_k\left[\erf\left(\frac{1-x_{ik}}{\sqrt{2}\ell_k}\right) - \erf\left(-\frac{x_{ik}}{\sqrt{2}\ell_k}\right)\right]\right\},
    \end{array}
\end{equation}
with $\erf$ being the error function, and $x_{ik}$ the $k$-th component of the observed input $\bx_{i}$.
The variance $\sigma_1^2$ is given by:
\begin{equation}
    \label{eqn:sigma_1}
        \begin{array}{ccc}
        \sigma_1^2 &:=& \E[\Q^{2}|\btheta, D_{n}] - (\E[\Q|\btheta, D_{n}])^{2} \\
        &=& \E[(\int_{\mathcal{X}}f(\bx)p(\bx)d\bx)^{2}||\btheta, D_{n}] - {\mu_{1}}^{2} \\
        &=& \E[\int_{\mathcal{X}}f(\bx)p(\bx)d\bx \int_{\mathcal{X}}f(\bx')p(\bx')d\bx'||\btheta, D_{n}] - {\mu_{1}}^{2} \\
        &=& \int_{\mathcal{X}}\int_{\mathcal{X}}\E[f(\bx)f(\bx')|\btheta, D_{n}]p(\bx)p(\bx')d\bx d\bx' - {\mu_{1}}^{2} \\
        &=& \int_{\mathcal{X}}\int_{\mathcal{X}}[k_n(\bx, \bx') + m_{n}(\bx)m_{n}(\bx')]p(\bx)p(\bx')d\bx d\bx' - {\mu_{1}}^{2}\\
        &=& \int_{\mathcal{X}}\int_{\mathcal{X}}k_n(\bx, \bx')p(\bx)p(\bx')d\bx d\bx' \\
        &=& \sigma_0^2 - \bm{\epsilon}_n^T\left(\mathbf{K}_n+\sigma^2\right)^{-1}\bm{\epsilon}_n,
        \end{array}
\end{equation}
where
\begin{equation}
    \label{eqn:sigma_0}
    \begin{array}{ccc}
    \sigma_0^2 &=& \int_{\mathcal{X}}\int_{\mathcal{X}} k(\bx,\bx')p(\bx)p(\bx')d\bx d\bx'\\
    &=& s^2 \prod_{k=1}^{d}(2\ell_{k}^{2}\sqrt{\pi})\Bigg\{\frac{-1}{\sqrt{\pi}} + \frac{1}{\sqrt{\pi}}\exp\left(\frac{-1}{2\ell_{k}^{2}}\right) + \\
    && \frac{1}{\sqrt{2}\ell_{k}}\erf\left(\frac{1}{\sqrt{2}\ell_k}\right)\Bigg\}.
    \end{array}
\end{equation}
\subsubsection{Quantification of the hypothetical state of knowledge about QoI}
\label{sec:hypothetical_Q}
To derive an analytical approximation of our hypothetical state of knowledge about the QoI, i.e., $p(\Q|\btheta,\mathbf{D}_n,\tilde{\bx}, \tilde{y})$, we proceed as in \sref{current_state_of_knowledge}, but 
with the remark that the posterior GP after adding the hypothetical observation will have mean function:
\begin{equation}
    \label{eqn:post_mean_hypothetical}
    \tilde{\mu}_{n+1}(\bx) = \mu_n(\bx) + k_n(\bx,\tilde{\bx})\frac{\tilde{y}-\mu_n(\tilde{\bx})}{k_n(\tilde{\bx},\tilde{\bx}) + \sigma^2},
\end{equation}
and covariance function:
\begin{equation}
    \label{eqn:post_cov_hypothetical}
    \tilde{k}_{n+1}(\bx,\bx') = k_n(\bx, \bx') - \frac{k_n(\bx,\tilde{\bx})k_n(\tilde{\bx},\bx')}{k_n(\tilde{\bx},\tilde{\bx}) + \sigma^2}.
\end{equation}
We get,
\begin{equation}
    \label{eqn:dist_2}
    p(\Q|\btheta, \mathbf{D}_{n},\tilde{\bx}, \tilde{y}) = \mathcal{N}(\Q|\mu_{2}(\tilde{\bx},\tilde{y}), \sigma^{2}_{2}(\tilde{\bx})).
\end{equation}

The mean $\mu_2(\tilde{\bx},\tilde{y})$ is:
\begin{equation}
    \label{eqn:mu_2}
    \begin{array}{ccc}
    \mu_{2}(\tilde{\bx},\tilde{y}) &:=& \E[\Q|\btheta, D_{n}, \tilde{\bx}, \tilde{y}] \\
    &=& \int_{\mathcal{X}}\tilde{\mu}_{n+1}(\bx)d\bx \\  
    &=& \mu_1 + \frac{\nu(\tilde{\bx})}{k_n(\tilde{\bx},\tilde{\bx}) + \sigma^2}\left(\tilde{y}-\mu_n(\tilde{\bx})\right)
    \end{array}
\end{equation}
with
\begin{equation}
    \label{eqn:nu}
    \nu(\tilde{\bx}) := \epsilon(\tilde{\bx})-\bm{\epsilon}_n^T\left(\mathbf{K}_n+\sigma^2\right)^{-1}\mathbf{k}_n(\tilde{\bx}),
\end{equation}
where $\epsilon(\tilde{\bx})$ as in \qref{epsilon_n} but with $\bx_i$ replaced by $\tilde{\bx}$.
Using the expression for the posterior covariance from \qref{post_cov_hypothetical} one can simplify $\sigma_2^2(\tilde{\bx})$ similar to the derivation in \qref{sigma_1} to get:
\begin{equation}
    \label{eqn:sigma_2}
    \begin{array}{ccc}
    \sigma_{2}^{2}(\tilde{\bx}) &:=& \E[\Q^{2}|\btheta, D_{n}, \tilde{\bx}, \tilde{y}] - (\E[\Q|\btheta, D_{n},\tilde{\bx},\tilde{y})^{2}\\
    &=& \int_{\mathcal{X}}\int_{\mathcal{X}}\tilde{k}_{n+1}(\bx,\bx')p(\bx)p(\bx')d\bx d\bx'\\
    &=&\sigma_{1}^{2} - \frac{\nu^2(\tilde{\bx})}{k_n(\tilde{\bx},\tilde{\bx})+\sigma^2}.
    \end{array}
\end{equation}

\subsubsection{Quantification of the expected information gain about the QoI}
\label{sec:information_gain}
Since both \qref{dist_1} and \qref{dist_2} are Gaussian, the KL divergence between the hypothetical and the current state of knowledge about the QoI conditional on the hyper-parameters, $G(\bx,\tilde{y};\btheta)$ of \qref{kld}, is analytically tractable~\cite{duchi2007derivations}, i.e.,
\begin{equation}
    \label{eqn:kld_Gaussian}
    \begin{array}{ccc}
    G(\bx,\tilde{y};\btheta) &=&\log\left(\frac{\sigma_{1}}{\sigma_{2}(\tilde{\bx})}\right) + \frac{{\sigma^{2}_{2}(\tilde{\bx})}}{2{\sigma^{2}_{1}}}  + \frac{{(\mu_{2}(\tilde{\bx},\tilde{y}) - \mu_{1})}^{2}}{2{\sigma^{2}_{1}}} - \frac{1}{2}.
    \end{array}
\end{equation}
Furthermore, $G(\bx,\tilde{y};\btheta)$ is a quadratic function of $\tilde{y}$, and $p(\tilde{y}|\tilde{\bx},\btheta,\mathbf{D}_n)$ is Gaussian, see \qref{point_predictive}.
Thus, we can analytically integrate $\tilde{y}$ out to obtain:
\begin{equation}
    \begin{array}{ccc}
    G(\tilde{\bx};\btheta) &=& \int_{-\infty}^\infty G(\tilde{\bx},\tilde{y};\btheta) p(\tilde{y}|\tilde{\bx},\btheta,\mathbf{D}_n)d\tilde{y}\\
    &=& \log\left(\frac{\sigma_1}{\sigma_{2}(\tilde{\bx})}\right) + \frac{1}{2}\frac{{\sigma_{2}}^{2}(\tilde{\bx})}{{\sigma_{1}}^{2}} - \frac{1}{2}\\
    && + \frac{1}{2}\frac{v(\tilde{\bx})^{2}}{\sigma_{1}^{2}(\sigma_{n}^{2}(\tilde{\bx}) + \sigma^{2})},
    \end{array}
\end{equation}
Finally, we take the expectation of $G(\tilde{\bx};\btheta)$ over the posterior of the hyperparameters, $p(\btheta|\mathbf{D}_n)$ of \qref{hyper_post}, using the MCMC samples $\left\{\btheta^{(s)}\right\}_{s=1}^S$ collected with the procedure described in~\cite{goodman2010ensemble,foreman2013emcee}. 
This yields: 
\begin{equation}
    \label{eqn:ekld_avg}
    \begin{array}{ccc}
    G(\tilde{\bx}) &=& \int G(\tilde{\bx};\btheta) p(\btheta|\mathbf{D}_n) d\btheta \\
    &\approx& \frac{1}{S}\sum_{s=1}^S G\left(\tilde{\bx};\btheta^{(s)}\right).
    \end{array}
\end{equation}

\subsubsection{Maximizing the expected information gain about the QoI}
At each stage of our BODE algorithm, we optimize the EKLD $G(\tilde{\bx})$ using Baysian global optimization (BGO) based on the augmented expected improvement (AEI)~\cite{huang2006global}.
This choice takes into account the noisy nature of the approximation of \qref{ekld_avg}, and it reduces the computational time compared to a brute force or a mutlistart-and-gradient-based-optimization approach.
See Algorithm~\ref{alg:bgo} for pseudocode.
In all our experiments, irrespective of the dimensionality, we use $T_n=20$ BGO iterations to optimize the EKLD.
\begin{algorithm}[htb]
\caption{Optimize the EKLD using BGO with AEI.}
\begin{algorithmic}[1]
\Require Initial number of EKLD evaluations $T_i$;
         maximum number of EKLD evaluations $T_n$;
         number of candidate designs $n_d$ for BGO;
         MCMC samples from the posterior of the hyperparameters $\left\{\btheta^{(s)}\right\}_{s=1}^S$;
         stopping tolerance $\gamma_i>0$.
    \State Evaluate $G(\tilde{\bx})$ using \qref{ekld_avg} at $T_i$ random points to generate training data, $\tilde{\bX}_{T_i} = \left\{\tilde{\bx}_1,\dots,\tilde{\bx}_{T_i}\right\}$ and $\mathbf{G}_{T_i} = \left\{\tilde{G}_1 = G(\bx_1),\dots,\tilde{G}_{T_i} = G(\bx_{T_i})\right\}$, for BGO.
    \State $t \leftarrow t_i$.
            \While {$t < T_n$} 
                \State Fit a standard GP on the input-output pairs $\tilde{\bX}_{t}$-$\tilde{\mathbf{G}}_t$ using maximum likelihood to approximate $G(\tilde{\bx})$.
                \State Generate a set of candidate test points $\hat{\bX}_{n_d}=\left\{\hat{\bx}_1,\dots,\hat{\bx}_{n_d}\right\}$ using Latin Hypercube Sampling (LHS) ~\cite{mckay2000}.
                \State Compute the AEI of all of the candidate points in $\hat{\bX}_{n_d}$.
                \State Find the candidate point ${\hat{\bx}}_j$ that exhibits the maximum AEI.
                \If{If the maximum AEI is smaller than $\gamma_i$}
                    \State Break.
                \EndIf
                \State Use \qref{ekld_avg} to evaluate $G(\tilde{\bx})$ at $\hat{\bx}_j$ measuring $\hat{G}_j = G(\hat{\bx}_j)$.
                \State $\tilde{x}_{t + 1} \leftarrow \hat{x}_j$.
                \State $\tilde{G}_{t + 1} \leftarrow \hat{G}_j$.
                \State $\bX_{t+1} \leftarrow \tilde{\bX}_{t}\cup\{{\tilde{\bx}}_{t + 1}\}$.
                \State $\mathbf{G}_{t+1} \leftarrow \mathbf{G}_{t}\cup\{\tilde{G}_{t + 1}\}$.
                \State $t\leftarrow t + 1$.
    \EndWhile
    \State return $\underset{{\tilde{\bX}_{T_n}}}{\arg\max} {\tilde{\mathbf{G}}}_{T_n}$.
\end{algorithmic}
\label{alg:bgo}
\end{algorithm}

\subsubsection{Complete BODE framework}
\label{sec:complete_framework}
In Algorithm~\ref{alg:ekld}, we provide pseudocode implementation of the proposed BODE framework.
The algorithm stops when a predetermined number of experiments have been performed.
Alternatively, one could stop the algorithm when the expected information gain is below a threshold.
\begin{algorithm}[htb]
\caption{Bayesian optimal design of experiments maximizing the expected information gain about the statistical expectation of a physical response.}
\begin{algorithmic}[1]
\Require Initially observed inputs $\bX_{n_i}$;
         initially observed outputs $\bY_{n_i}$;
         maximum number of allowed experiments $N$.
    \State $n \leftarrow n_i$.
            \While {$n < N$}
                \State Sample from the posterior of the hyperparameters, \qref{hyper_post}, to obtain $\left\{\btheta^{(s)}\right\}_{s=1}^S$.
                \State Find the next experiment $\bx_{n+1}$ using Algorithm \ref{alg:bgo} to solve \qref{iaf_def}.
                \State Evaluate the objective at $\bx_{n+1}$ measuring $y_{n+1} = f(\bx_{n+1})$.
                \State $\bX_{n+1} \leftarrow \bX_{n}\cup\{\bx_{n+1}\}$.
                \State $\bY_{n+1} \leftarrow \bY_{n}\cup\{y_{n+1}\}$.
                \State $t\leftarrow t + 1$.
        \EndWhile
\end{algorithmic}
\label{alg:ekld}
\end{algorithm}

\section{Results}
\label{sec:results}
We apply the methodology on two one-dimensional mathematical 
functions (synthetic problems), a three-dimensional problem, and a five-dimensional problem. 
For the first two synthetic problems the input domain simply becomes $[0, 1]$ whereas for the third synthetic problem the input domain is $[-2, 6]^{3}$. 
The inputs for the five dimensional numerical example lie in the hyper-cube $[0, 1]^{5}$. 
The number of initial data points is denoted by $n_{i}$. 
The number of initial data points is taken as low as possible for the numerical examples. 
In most literature, as a rule of thumb, $10d$ number of initial samples are used. 
We resort to using lesser number of initial data points to test the performance of the methodology when it starts from the low-sample regime.
Readers interested in the problem of the optimal selection of initial data size can refer to the work of S\`obester et. al.~\cite{sobester2005design} where the authors discuss the problem in the context of optimization.
The problem of selecting an optimal number of initial points is beyond the scope of the work presented here. 

\subsection{Synthetic problem no. 1}
\label{sec:toy_1}
Consider the function
\begin{equation}
    f(x) = 4\left( {1 - \sin\left(6x + 8{e^{6x - 7}}\right)}\right),
    \label{eqn:toy_1}
\end{equation}
defined on $[0,1]$.
This function is smooth throughout its domain, but it exhibits two local minima.
We will apply our methodology to estimate the statistical expectation:
$$
Q[f] = \int_0^1 f(x)dx.
$$
The true value of $Q[f]$ is analytically available, $Q[f] = -1.3599$.
We apply our methodology to this problem starting from $n_i=3$ and sample a total of $N = 28$ points.
The number of MCMC chains for the results shown below is six, and the number of steps per chain is 500.
For further details on the MCMC part of training the GP, we refer the readers to ~\cite{foreman2013emcee,goodman2010ensemble}. 

\begin{figure}[!htbp]
\centering
    \subfigure[]{
        \includegraphics[width=1.\columnwidth]{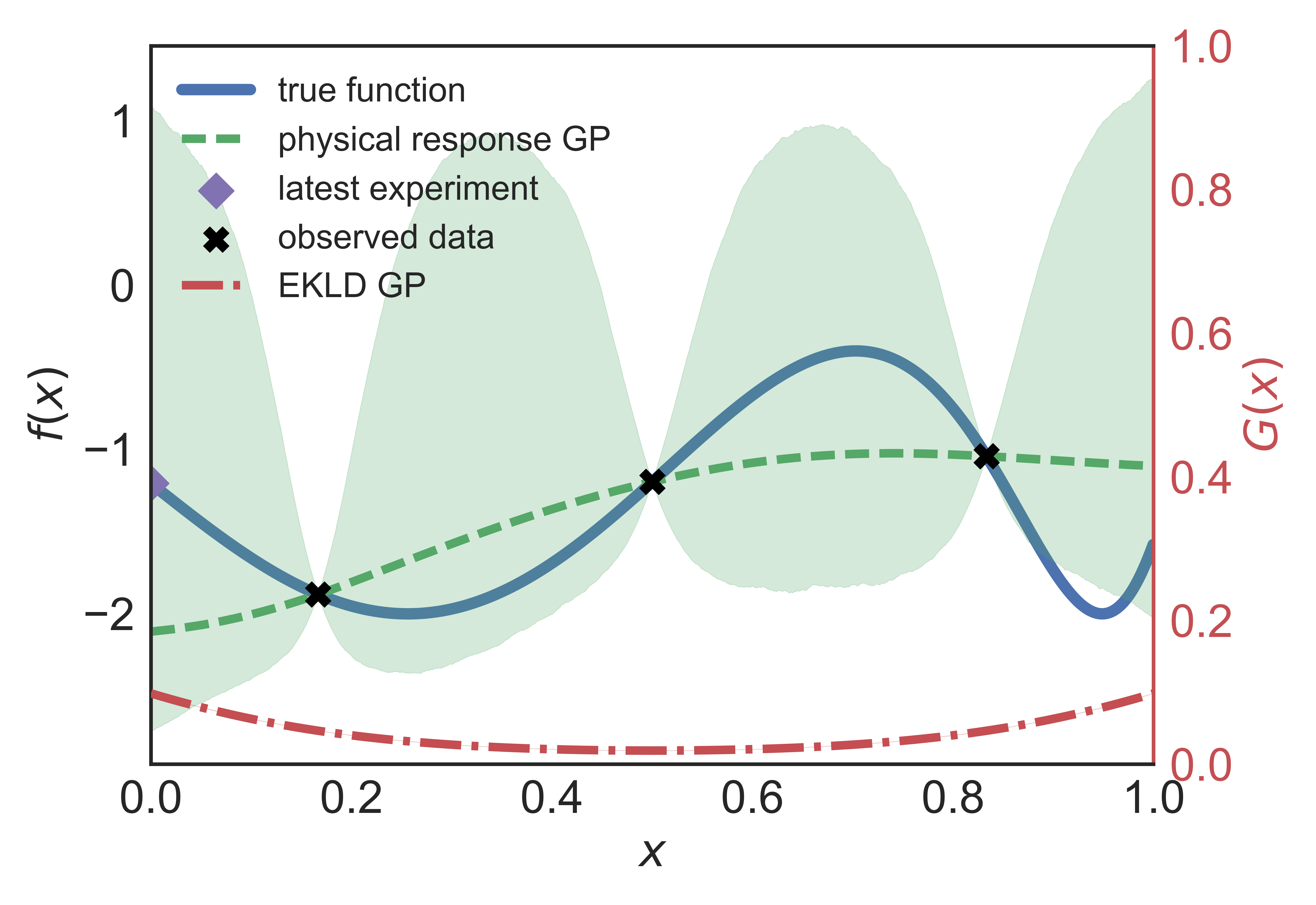}
    }
    \subfigure[]{
        \includegraphics[width=1.\columnwidth]{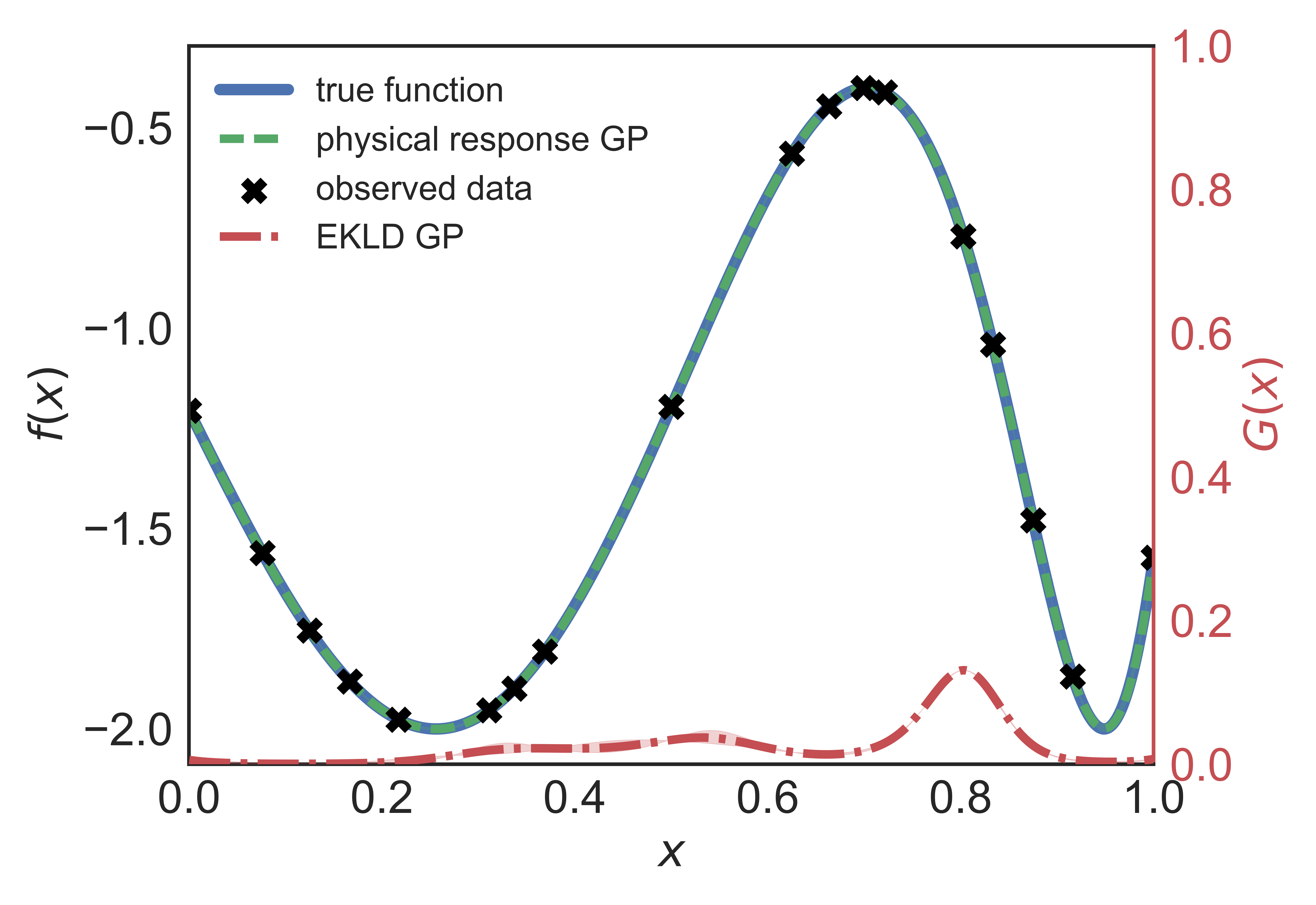}
    }
    \caption{One-dimensional synthetic problem ($n_{i}=3$).
        Subfigures ~(a) and (b) show the state of the function (1st iteration) at the start
        and the end (15th iteration) of the algorithm.
    }
    \label{fig:toy_1_unif}
\end{figure}

Figs.~\ref{fig:toy_1_unif}~(a) and~(b) show the initial and final state of Algorithm~\ref{alg:ekld}.
The thick blue line represents the true function $f$, \qref{toy_1}.
The black crosses are the observed data at the given stage.
In subfigure~(a), the next experiment selected by maximizing the EKLD, see Algorithm \ref{alg:bgo}, is corresponds to the purple diamond.
The mean of the GP fit to the expected information gain $G(\tilde{\bx})$ constructed by BGO in Algorithm~\ref{alg:bgo}.
The predictive mean of the EKLD is shown by the dotted light blue line.
This dotted line represents the response surface of the EKLD after the BGO has ended and the red shaded area around it represents the uncertainty (2.5 percentile and  97.5 percentile) around it.
As expected, the mean of the EKLD is very small or close to zero at points where experiments have been performed.
Thus, the point selected by the methodology (purple diamond) is located in the input space where the EKLD has high mean.
The posterior mean of the GP of the black-box function is represented by the dashed bottle-green line.
The bottle-green shaded area represents the uncertainty (2.5 percentile and  97.5 percentile) around it.
The final set of inputs, space-filling, selected by the methodology can be seen in \fref{toy_1_unif}~(b). 
\fref{toy_box}~(a) shows the $p(\Q|\mathbf{D}_n)$ plotted against the number of data samples while showing 
convergence towards the true value of $\Q[f]$. The gradual reduction of predictive uncertainty of $\Q$ from the initial to the final stage of the algorithm is seen in \fref{toy_box}~(a).

\subsection{Synthetic problem no. 2}
\label{sec:toy_2}
We consider the following Gaussian mixture function to test and validate our methodology further.
\begin{equation}
    \label{eqn:toy_2}
    \begin{array}{ccc}
    f(x) &=&  \frac{1}{\sqrt{2\pi}s_{1}}\exp\left\{-\frac{(x-m_{1})^{2}}{2{s_{1}}^2}\right\}\\
    && + \frac{1}{\sqrt{2\pi}s_{2}}\exp\left\{-\frac{(x-m_{2})^{2}}{2{s_{2}}^2}\right\},
    \end{array}
\end{equation}
where $m_{1}=0.2$ and $s_{1}=0.05$, $m_{2}=0.8$ and $s_{2}=0.05$.
As can be seen from \qref{toy_2}, the function is a sum of probability densities of two Gaussian distributions. 
The notoriety of the function lies in two relatively sharp but smaller areas of high magnitude.
The true value of $\Q[f]$ is analytically available, $\Q[f] = 2.0$.
We apply our methodology to this problem starting from $n_i=3$ and sample another 25 points.
The final state of sampling can be seen in \fref{toy_2_unif}~(b), which shows a fairly equally spaced spread of designs.
It is important to note that \fref{toy_2_unif}~(b) can mislead the reader into perceiving the sampling to be less dense in the areas where the function is sharply peaked.
This is an illusion due to the starkly varying ordinates of the sampled points near the peaks of the function.
The convergence of the estimated mean to the true value of $\Q[f]$ and the reduction in uncertainty around the $\Q[f]$ can be seen in \fref{toy_box}~(b).

\begin{figure}[!htbp]
\centering
    \subfigure[]{
        \includegraphics[width=1.\columnwidth]{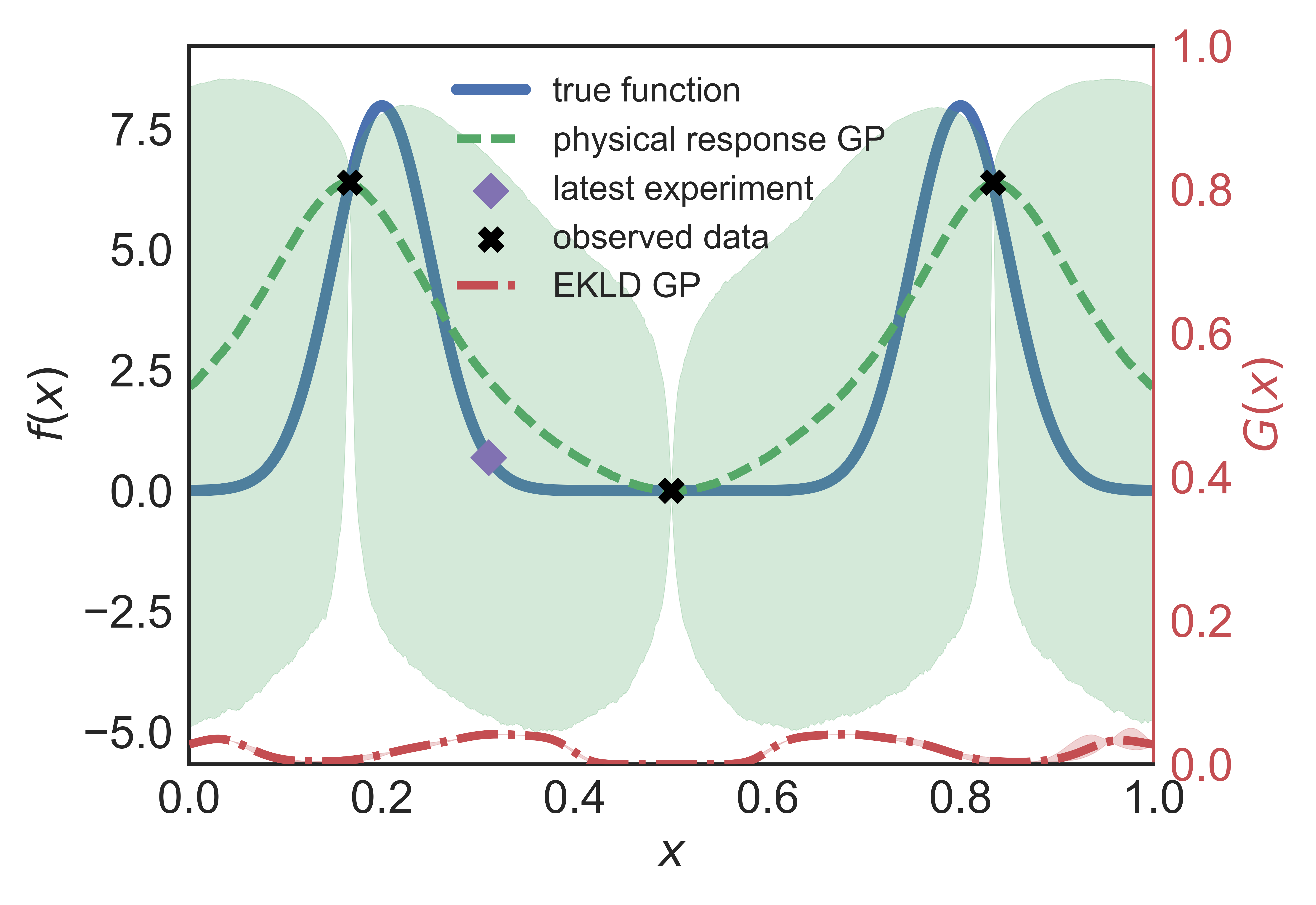}
    }
    \subfigure[]{
        \includegraphics[width=1.\columnwidth]{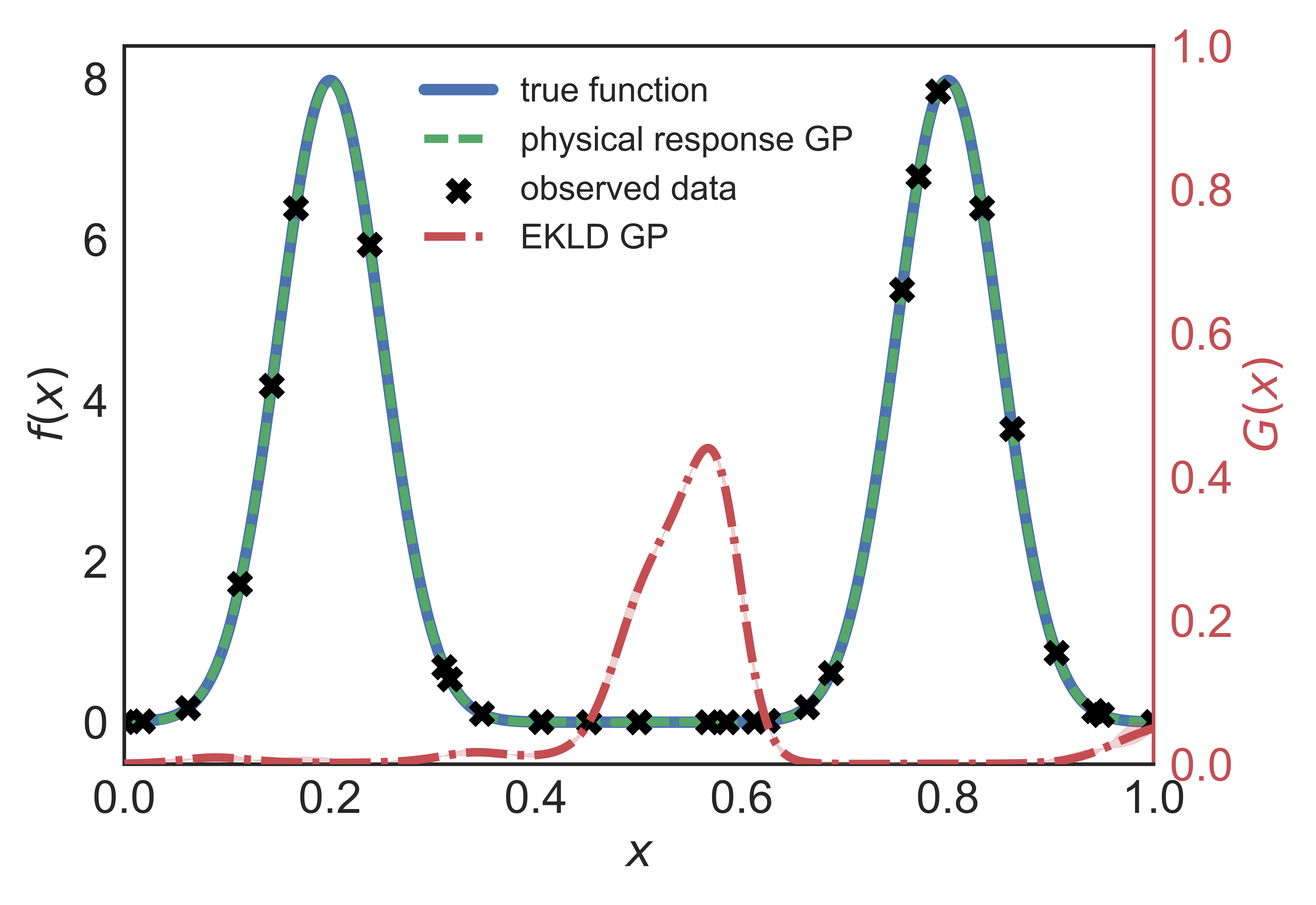}
    }
    \caption{One-dimensional synthetic example ($n_{i}=3$).
        Subfigures ~(a) and (b) show the state of the function at the start (1st iteration)
        and the end (25th iteration) of the algorithm.
    }
    \label{fig:toy_2_unif}
\end{figure}
In \fref{toy_2_unif} (a) and (b), the EKLD is shown by the dotted line and the true function is shown by the dashed red line. The solid line represents the mean of the GP model and the orange shaded areas around it represent the 2.5th and the 97.5th percentile of the GP.
We plot the relative maximum mean EKLD as a function of the number of samples in \fref{toy_ekld} for both the synthetic functions. This relative maximum EKLD is the ratio of the maximum predictive mean of the EKLD 
for the current iteration and the overall maximum predictive mean of the EKLD obtained across all iterations.
The plots in \fref{toy_ekld} show a characteristic typical of BODE functions i.e. of increasing in magnitude for 
the first few iterations and then falling sharply. 
This predicted mean value of the EKLD asymptotically goes to zero for both the synthetic 
functions here. The number of MCMC chains for the results shown below is six, and the number of steps per chain is 500.
\begin{figure}[!htbp]
\centering
    \subfigure[]{
        \includegraphics[width=1\columnwidth]{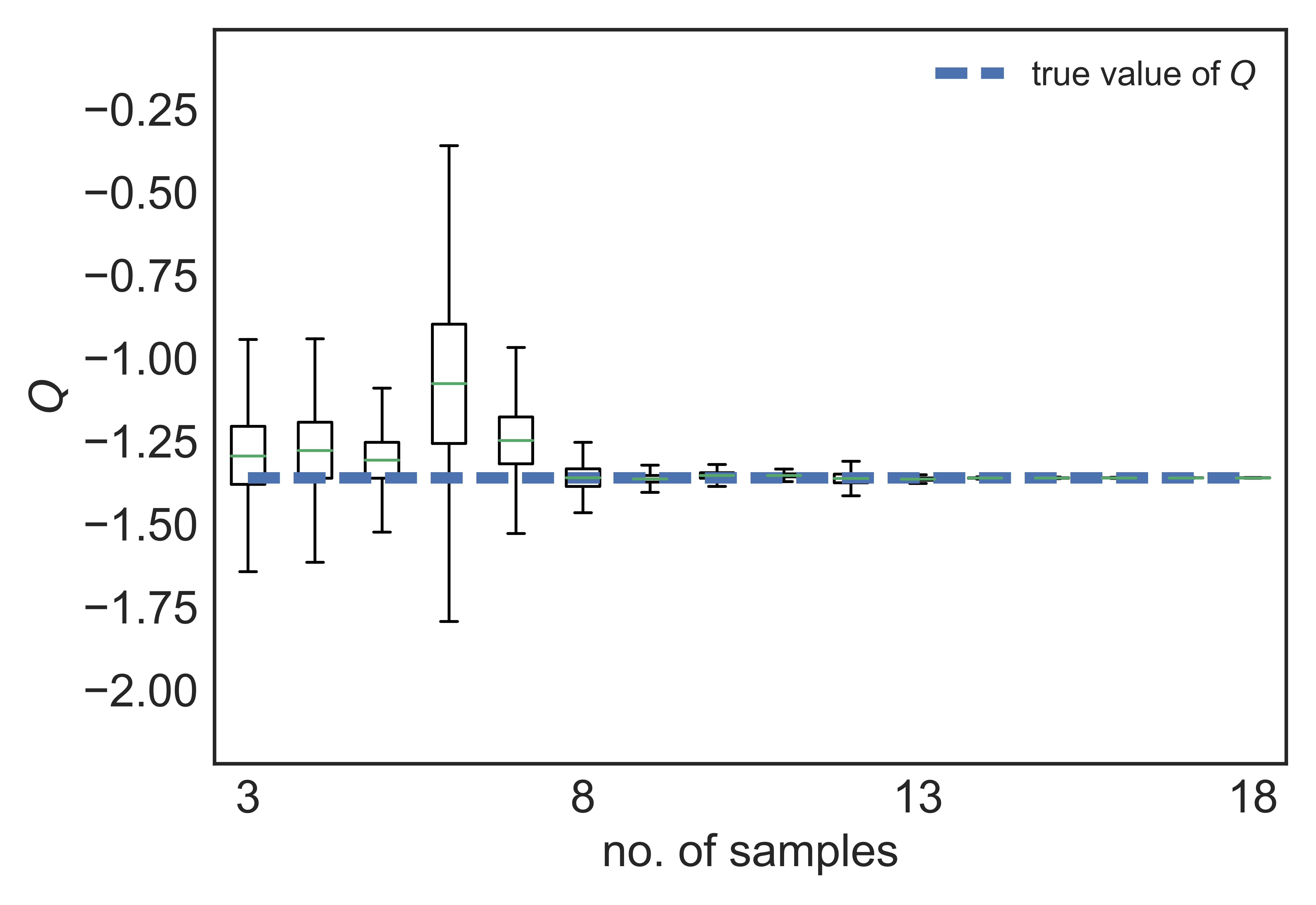}
    }
    \subfigure[]{
        \includegraphics[width=1\columnwidth]{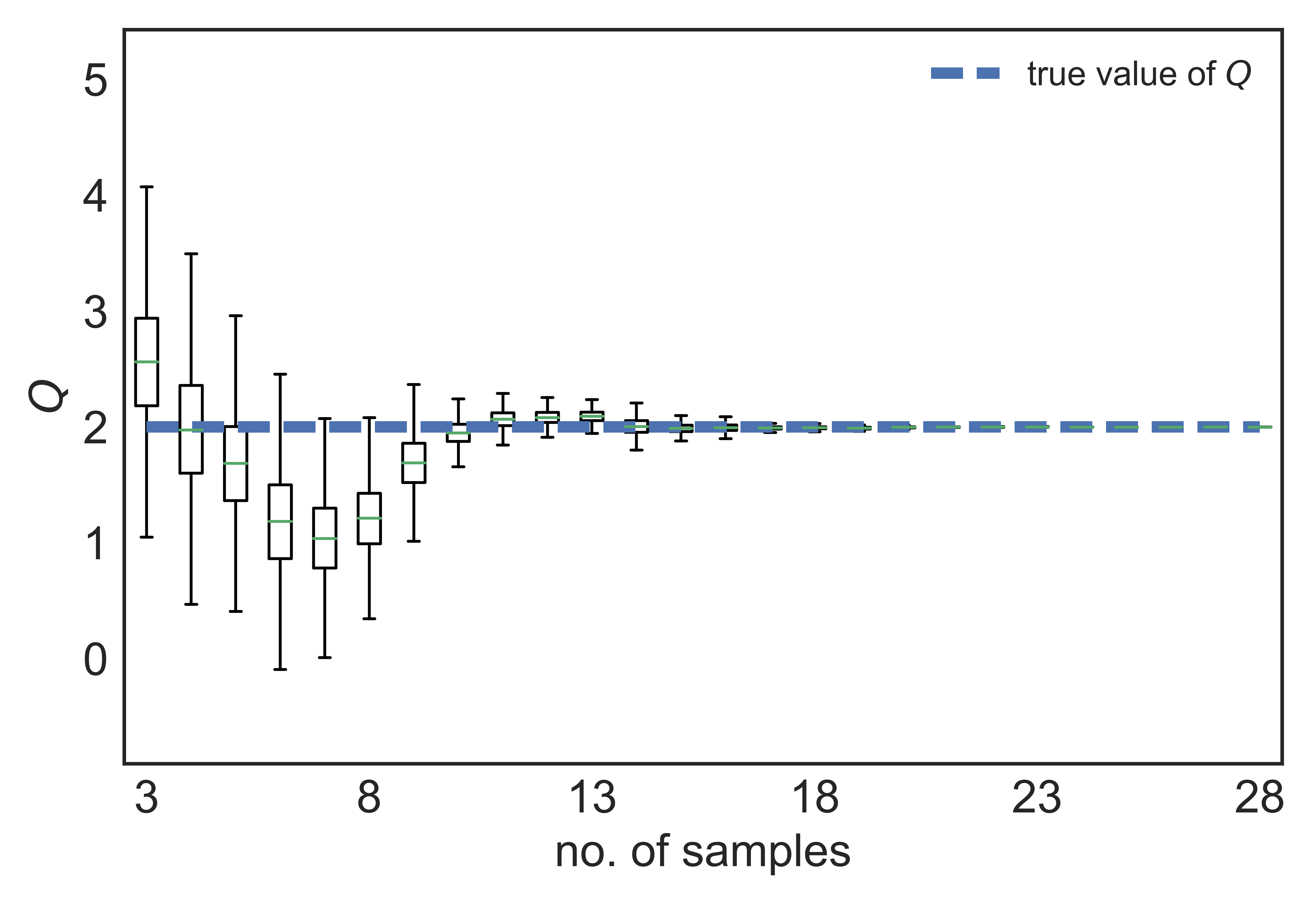}
    }
    \caption{One-dimensional synthetic examples.
        Subfigures ~(a) and (b) show the convergence to the true expectation of the function and the reduction
        in uncertainty about the QoI after the end of the algorithm, 
        for synthetic problem no. 1 ($n_{i}=3$) and synthetic problem no. 2 ($n_{i}=3$) respectively.
    }
    \label{fig:toy_box}
\end{figure}

\begin{figure}[!htbp]
\centering
    \subfigure[]{
        \includegraphics[width=1\columnwidth]{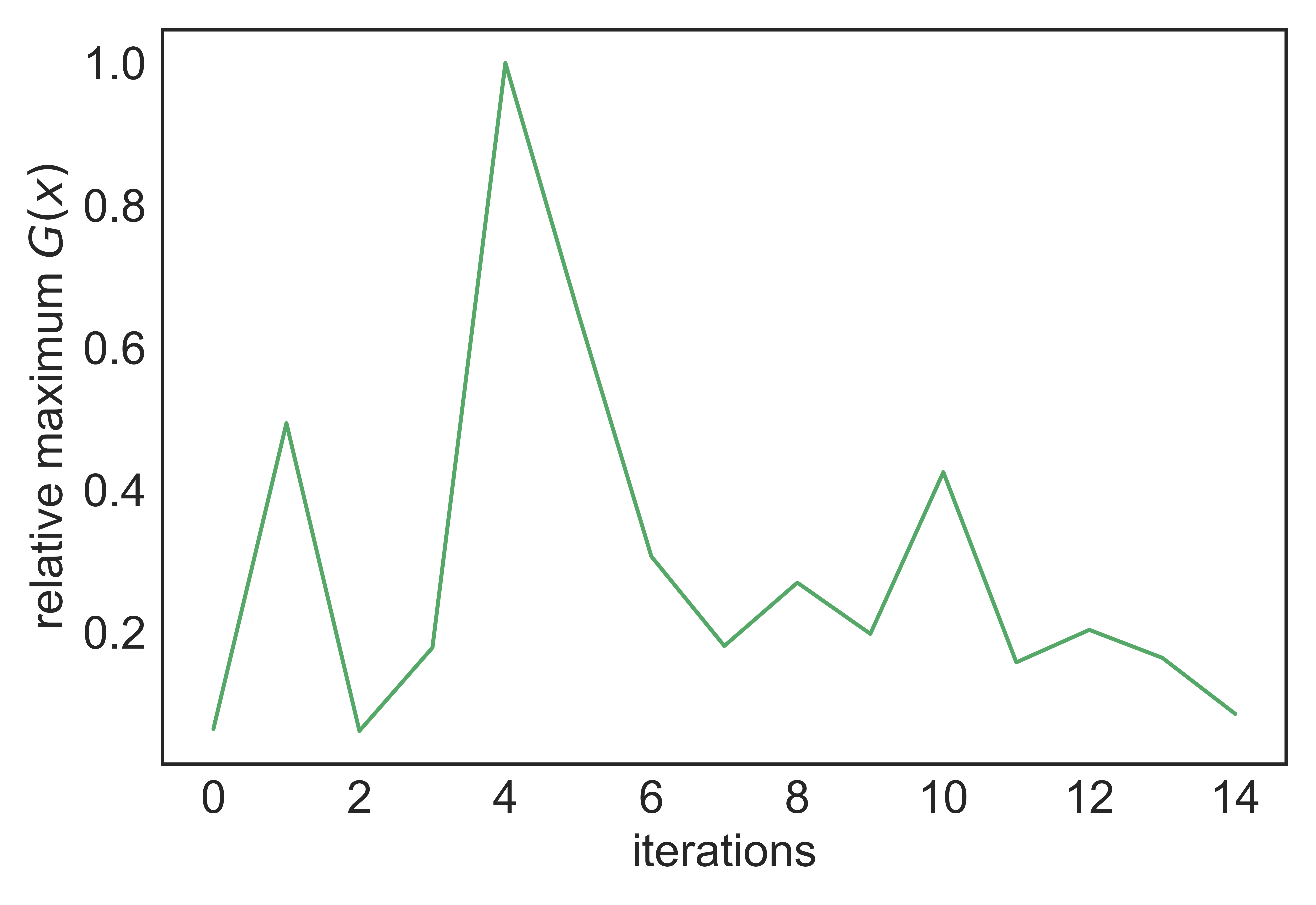}
    }
    \subfigure[]{
        \includegraphics[width=1\columnwidth]{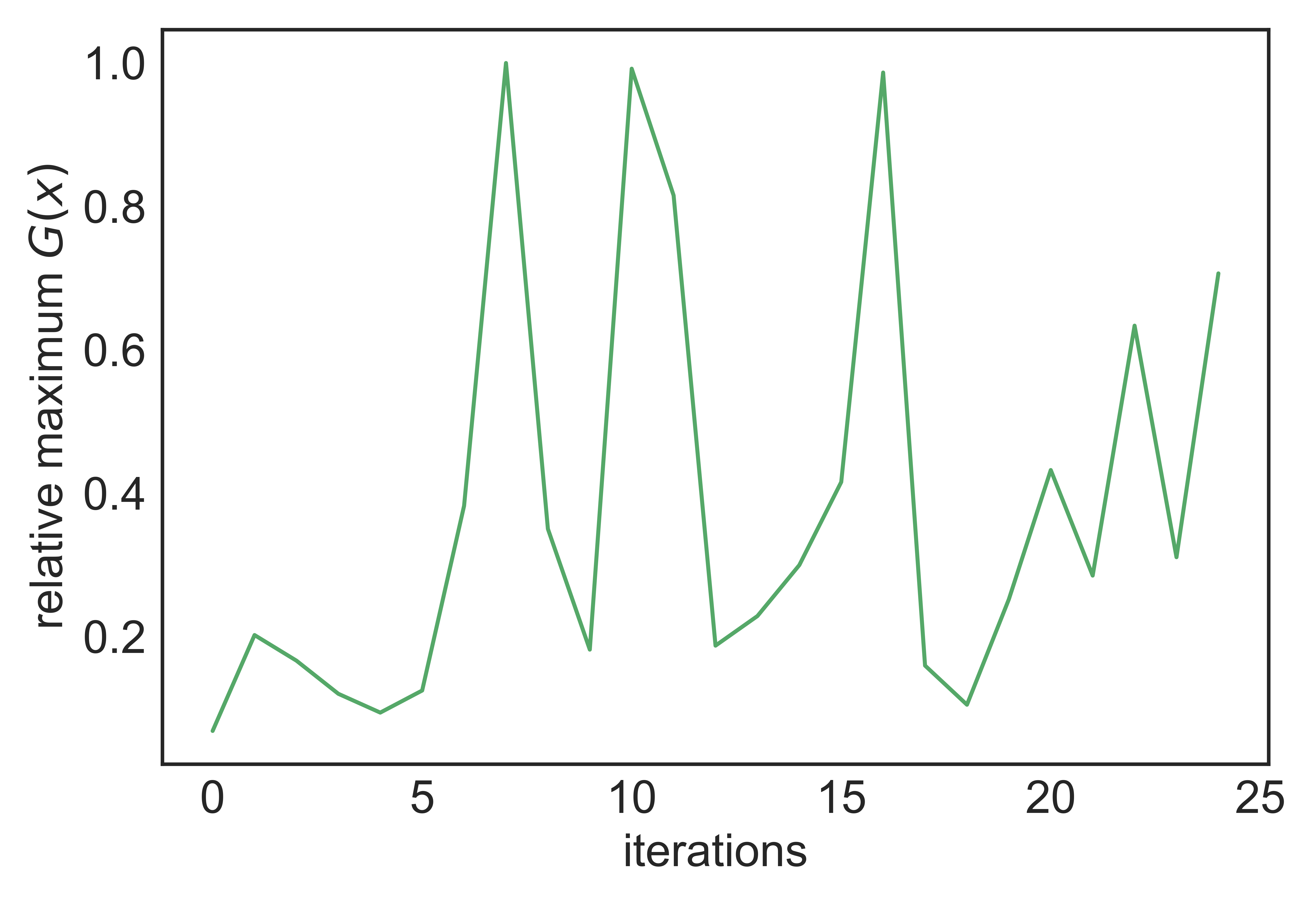}
    }
    \caption{One-dimensional synthetic examples.
        Subfigures ~(a) and (b) show the predictive mean of the EKLD, 
        for synthetic problem no. 1 ($n_{i}=3$) and synthetic problem no. 2 ($n_{i}=4$) respectively.
    }
    \label{fig:toy_ekld}
\end{figure}

\subsection{Synthetic problem no. 3}
\label{sec:toy_3}
We consider the following three dimensional function from \cite{dette2010generalized} to test and validate our methodology further. 
\begin{eqnarray}
    f(\bx)&=&4(x_{1} + 8x_{2} - 8x_{2}^{2} - 2)^{2} + (3-4x_{2})^{2} \nonumber \\
    &&+ 16\sqrt{x_{3}+1}(2x_{3}-1)^{2}.
\label{eqn:toy_3}
\end{eqnarray}
The major difference between this function \qref{toy_3} and the the first two synthetic examples is the dimensionality of the problem.
The true value of $\Q[f]$ is analytically available, $\Q[f] = -0.7864$.
We apply our methodology to this problem starting from $n_i=2$ and sample another 30 points.
\fref{toy_3_unif}~(b) shows that the methodology started with a highly uncertain estimate of the true value and eventually converged to a sharp peaked Gaussian distribution around the true value. The approximation to $\Q[f]$ at each stage of the algorithm is shown in \fref{toy_3_unif}~(b). The gradual reduction in uncertainty around $\Q[f]$ also can be seen in \fref{toy_3_unif}~(b). 
\fref{toy_3_unif}~(a) demonstrates how the relative EKLD fluctuates while seemingly approaching zero.

\begin{figure}[!htbp]
\centering
    \subfigure[]{
        \includegraphics[width=1\columnwidth]{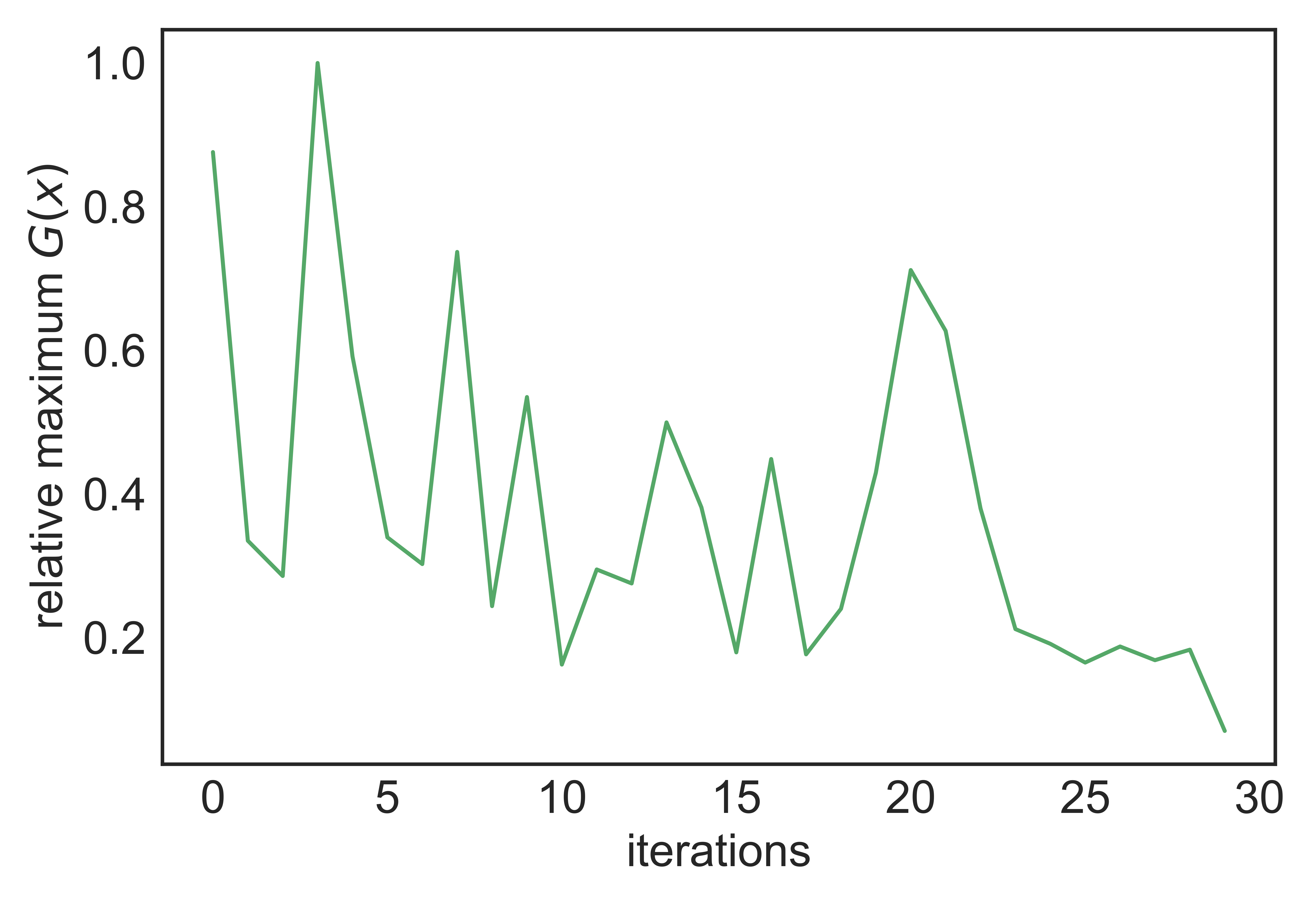}
    }
    \subfigure[]{
        \includegraphics[width=1\columnwidth]{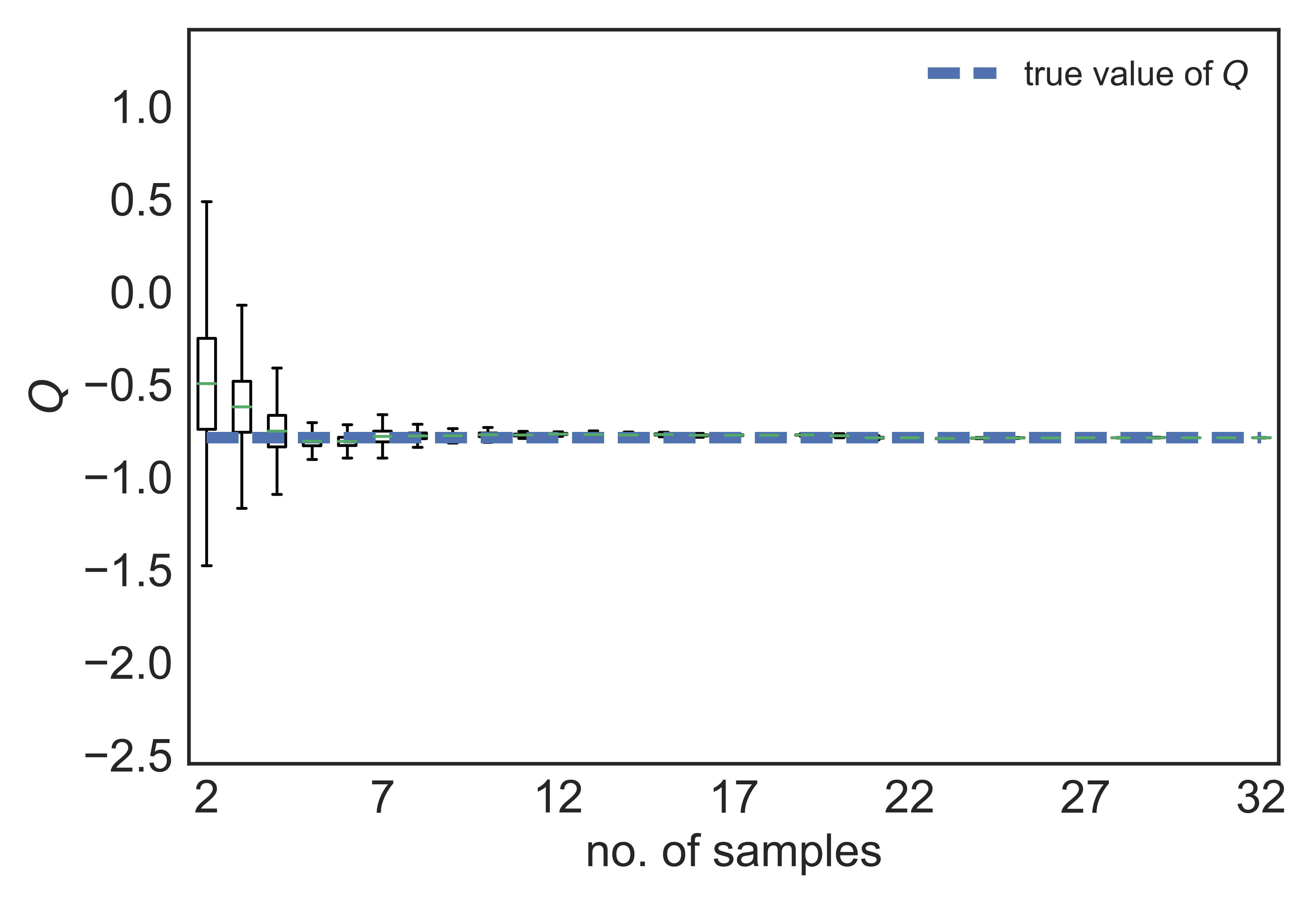}
    }
    \caption{Three-dimensional synthetic example ($n_{i}=2$).
        Subfigure ~(a) shows the decay of the EKLD from the 1st iteration
        to the end of the 30th iteration of the algorithm.
        Subfigures ~(b) show the convergence to the true value of the QoI
        respectively.
    }
    \label{fig:toy_3_unif}
\end{figure}

\subsection{Synthetic problem no. 4}
\label{sec:toy_4}
The following five dimensional function is taken from \cite{knowles2006parego}. 
\begin{eqnarray}
    f(\bx)&=& 10\sin(\pi x_{1}x_{2}) + 20 (x_{3}-5)^{2} +10x_{4} + 5x_{5}. 
\label{eqn:toy_4}
\end{eqnarray}
This function \qref{toy_4} is reasonably high-dimensional and challenging due to the non-linear input-output relation. 
The true value of $\Q[f]$ is analytically available, $\Q[f] = 0.3883$.
We apply our methodology to this problem starting from $n_i=20$ and sample another 45 points.
\fref{toy_4_unif} (a) demonstrates how the mean of the relative EKLD tends to approach zero by the end of the sampling process.  
The iteration-wise convergence of the $\Q[f]$ to its true value is shown in \fref{toy_4_unif} (b).
\fref{toy_4_unif} (b) can present an illusion to the reader as it shows that the mean of the QoI is very close to the true value at the start of sampling. This is misleading because of the relatively large variance around the mean which means that the methodology is not confident of being close to the true value. As a result of this it can be seen, in the subsequent iterations, that the mean of the QoI goes to either side of the true value with a gradual decrease in variance. This might happen due to the methodology discovering different modes of the underlying function. As more data are accumulated, the uncertainty around the estimate decreases.
\begin{figure}[!htbp]
\centering
    \subfigure[]{
        \includegraphics[width=1\columnwidth]{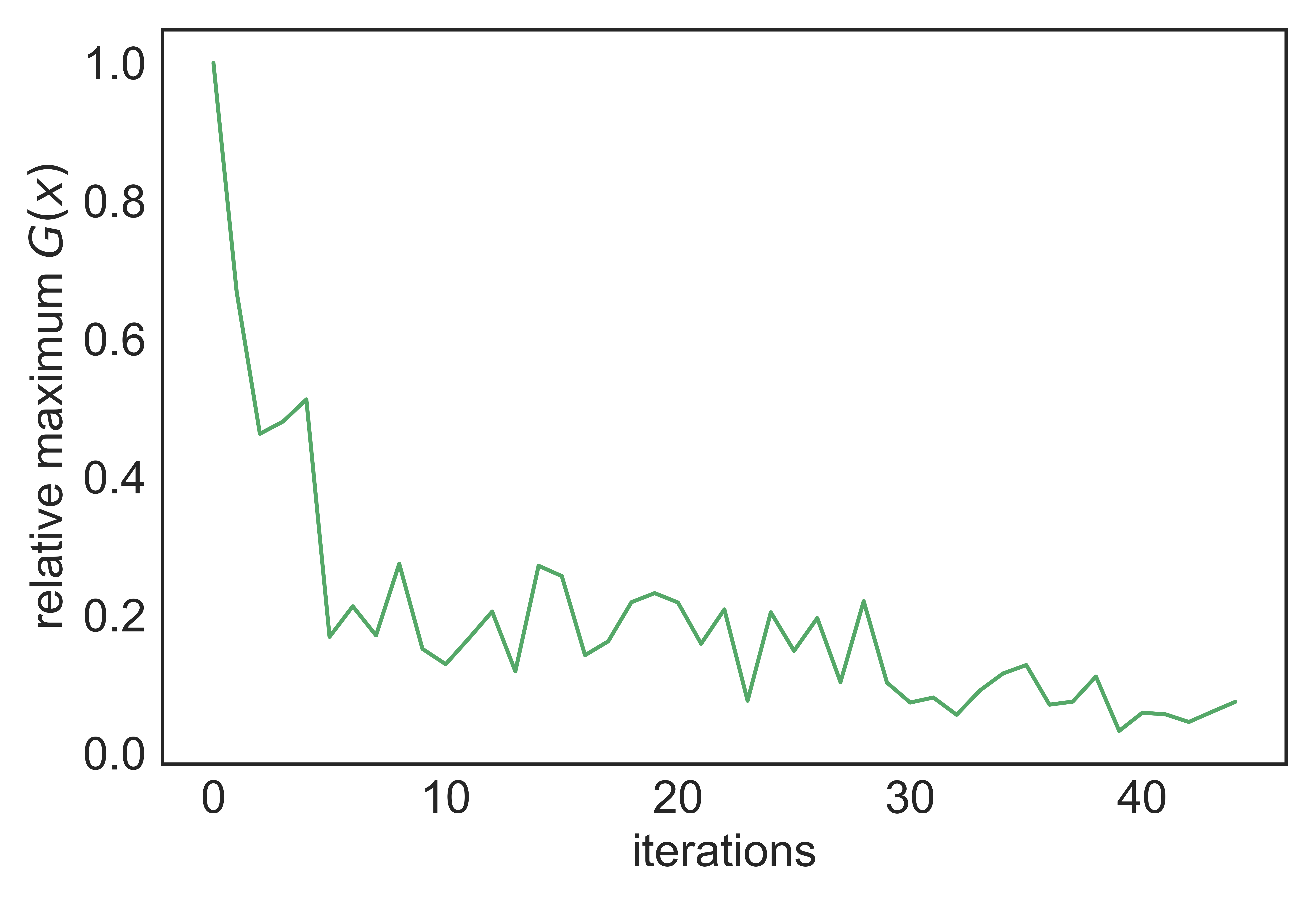}
    }
    \subfigure[]{
        \includegraphics[width=1\columnwidth]{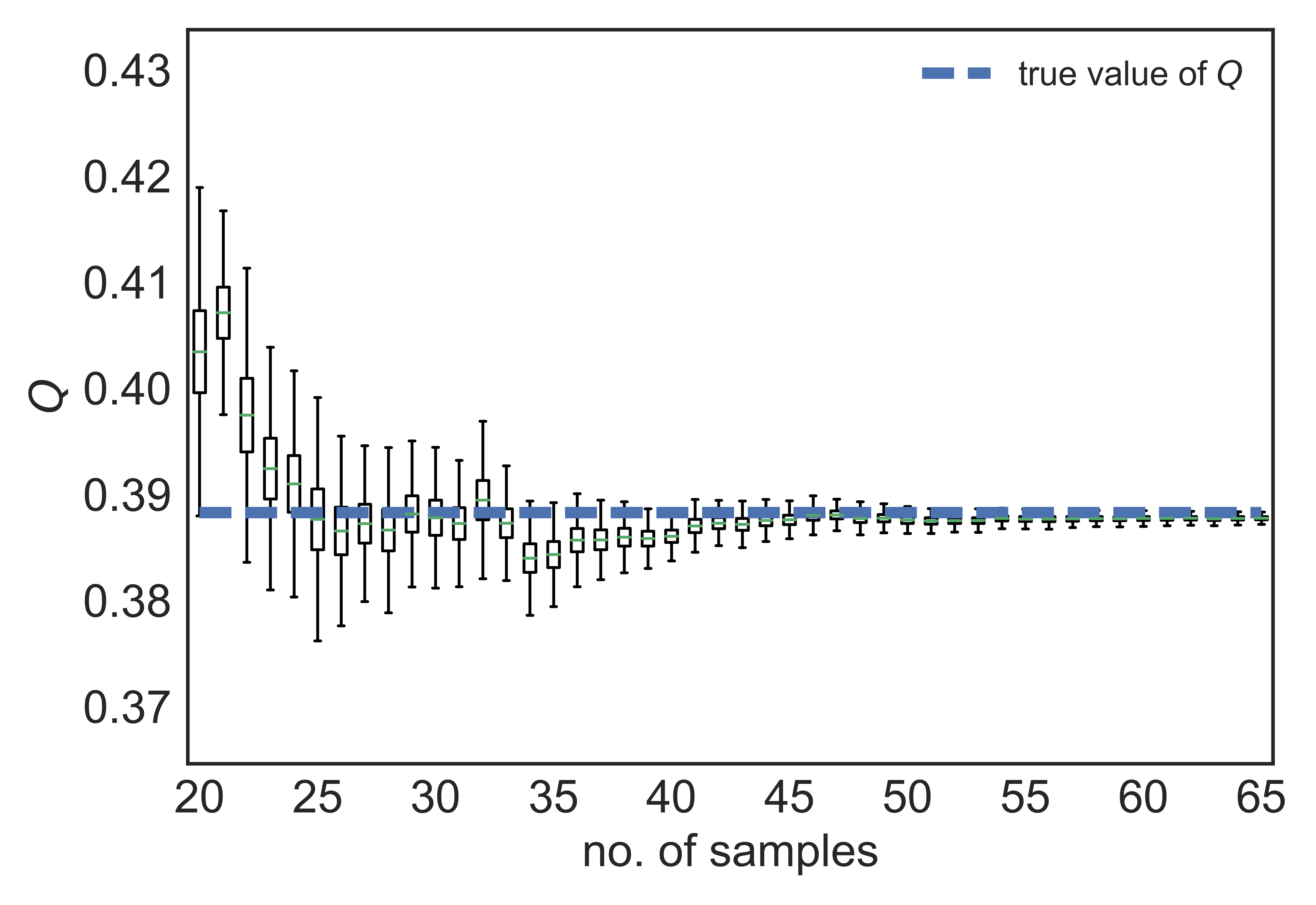}
    }
    \caption{Five-dimensional synthetic example ($n_{i}=20$).
        Subfigure ~(a) shows the decay of the EKLD from the 1st iteration
        to the end of the 45th iteration of the algorithm.
        Subfigure ~(b) shows convergence to the true value of the QoI.
    }
    \label{fig:toy_4_unif}
\end{figure}

\subsection{Steel wire drawing problem}
\label{sec:wire}
The wire drawing process aims to achieve a required reduction in the cross section of the incoming wire, while 
aiming to monitor or optimize the mechanical properties of the outgoing wire. 
The incoming wire is passed through a series of dies (8 dies) to achieve an overall reduction in wire diameter.
Each pass reduces the cross section of the incoming wire.
The authors refer the reader to Section 3.4 in \cite{pandita2018stochastic} to gain further information about the modeling of the wire drawing process using finite element method.
The wire drawing process here is represented by an expensive computer code of which only a small number of evaluations are possible. The frictional work per Tonne (FWT) is one of the outputs of the expensive code. The statistical expectation value of the FWT is of importance for various stakeholders as the work done by the friction on the passing wire determines the power consumed, the wear on the final wire, etc. The FWT is the aggregate of the frictional work done at each pass.
In our problem, we consider the die angle as design variables for each pass. The outgoing diameters at each pass are fixed to reasonable values. Thus, we deal with a total of 8 design variables.  
We start the methodology with 20 initial data points and add another 80 samples.
We approximate the true value of the expectation of FWT, by averaging the outputs at 6,000 designs generated by Latin-hypercube sampling (LHS), as $\Q[FWT] \approx 0.2694$.
The results in \fref{wire_unif} show the gradual convergence of the methodology's mean estimate of the QoI towards the approximated true value.
\fref{wire_unif} (a) shows the mean and variance of the expectation of FWT as the mean approaches the approximate true value while the variance around it decreases gradually. 
The reduction in variance around the QoI from the start of the sampling to the end can be seen in \fref{wire_unif} (b).
This is intuitive as the number of collected samples increases, the variance around the QoI decreases. 
The comparison of the performance of the EKLD to that of the US is seen in \fref{wire_unif} (c). The mean of the statistical expectation value of FWT for the EKLD converges to the approximate true value as more samples are added, while that for the US makes gradual drifts either side of the approximate true value. The US requires more samples to approach the approximate true value. This difference may be explained by the context specific functional form of the derived EKLD compared to the agnostic US which, although is a reduced form of the KLD in the design variables, seems to be slower in higher dimensions.
\begin{figure}[!htbp]
\centering
    \subfigure[]{
        \includegraphics[width=1\columnwidth]{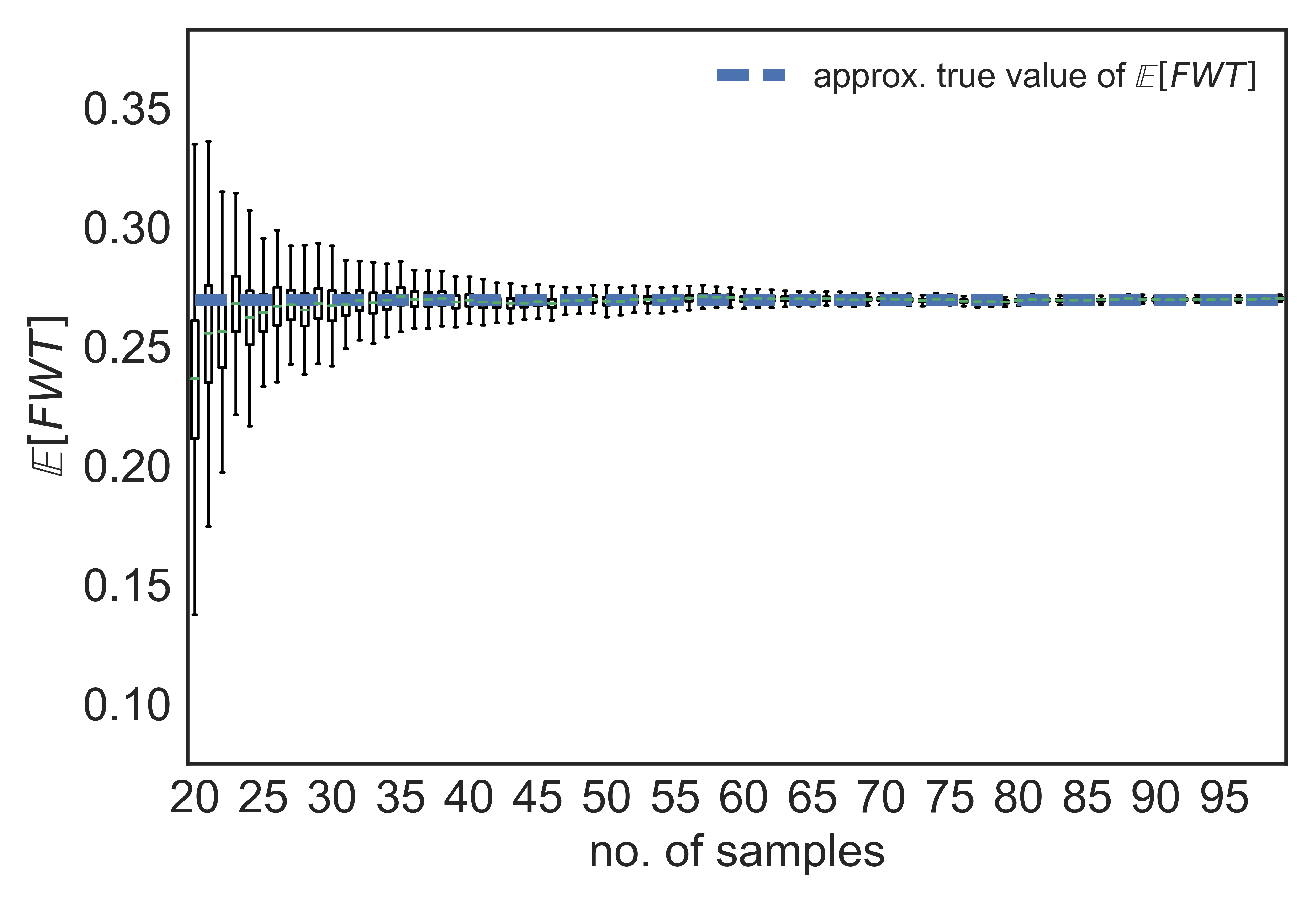}
    }
    \subfigure[]{
        \includegraphics[width=1\columnwidth]{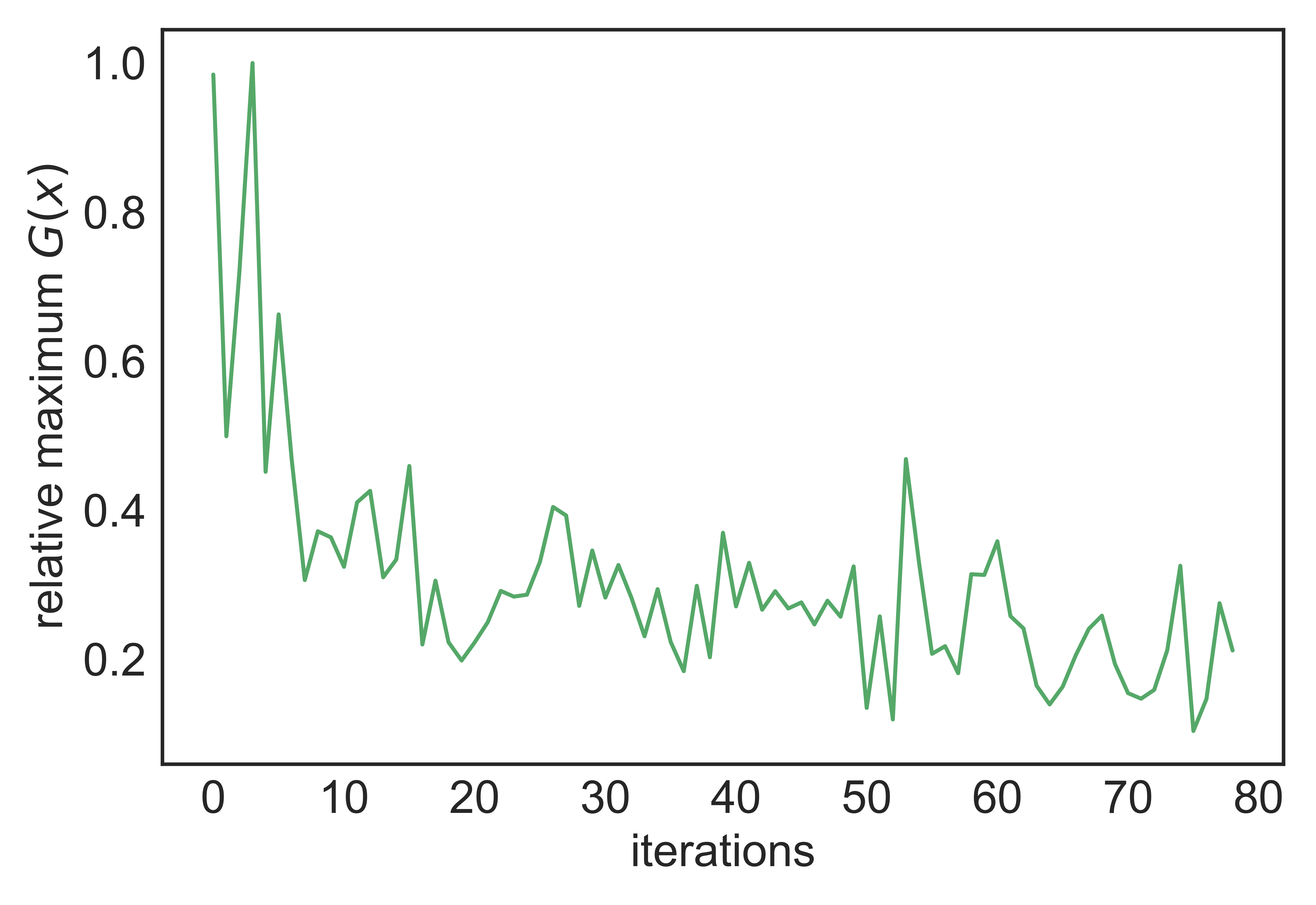}
    }
    \caption{Wire drawing problem ($n_{i}=20$) after 80 iterations. 
       }
    \label{fig:wire_unif}
\end{figure}

\subsection{Comparison with Uncertainty Sampling}
\label{sec:comp}
As a demonstration of the performance of the methodology in contrast to a ubiquitous state-of-the-art
sampling technique, namely uncertainty sampling (US), the methodology
is tested on the synthetic examples given in \sref{toy_1}, \sref{toy_2}, \sref{toy_3} and \sref{toy_4}.
The uncertainty sampling technique works on the principle of reducing the uncertainty
around the predictive response surface. Interestingly it has been shown that maximizing the information gain in the parameters reduces to uncertainty sampling under certain assumptions \cite{mackay1992information}. 
Moreover, US, in its functional form, as an IAF is agnostic to the context (QoI) in the problem. 
Hence, it serves as an ideal benchmark to compare with the EKLD. 
An explanation of the US methodology is as follows.
The methodology selects a design with the maximum magnitude of predictive variance and follows this procedure until it sequentially acquires the required number of samples. 
The surrogate modeling process for the US works the same way as for the EKLD. 
The overall algorithm remains the same as Algorithm \ref{alg:ekld}, but for the change in the sampling criterion. 

The convergence to the QoI for the synthetic problem in \sref{toy_1} and \sref{toy_2} is seen in \fref{toy_comp_1} (a) and \fref{toy_comp_1} (b) respectively.
Overall, the two methodologies converge to the true value within reasonable time of one another. 
With the two peaked one-dimensional function of \sref{toy_2}, the EKLD takes more iterations to converge as seen in \fref{toy_comp_1} (b).
The US can be seen as being quicker in reaching very close to the true value of the QoI compared to the EKLD for the synthetic problem no. 2 whereas EKLD takes slightly fewer iterations to estimate the true statistical expectation value for synthetic problem no. 1.

As the complexity of the problems increases, convergence for the EKLD becomes quicker compared to US as shown in \fref{toy_comp_1} (a), (b) and (c).
With the three-dimensional problem \fref{toy_comp_2} (a), the mean estimate of the QoI for the EKLD converges after 20 samples have been collected.
For the same problem, US takes almost 30 samples to converge.
This saving of almost 10 samples could be useful in engineering problems where each sample is collected at the expense of thousands of dollars of effort or a computational burden of multiple days.

For the five-dimensional synthetic problem, \fref{toy_comp_2} (b) shows how the EKLD starts to approach the true value of the QoI as the number of iterations increases, whereas US tends to shows jaggedness in its patterns of convergence.
After 65 samples have been collected US shows convergence, but convergence can be seen for the EKLD as early as the addition of the 45th sample.
This observation is further strengthened by looking at the decay of the EKLD in \fref{toy_4_unif} (b).
The comparison in \fref{toy_comp_2} (b) highlights the capability of the methodology to infer the QoI in a limited number of iterations.
This is useful in the context of problems with expensive black-box functions where each evaluation of the expensive function has a very high cost. 
Moving on to the wire-problem in \fref{wire_comp}, it can be seen that the convergence to the approximate true value is achieved by the EKLD and US albeit with more samples for US.

Another important feature of the comparisons in \fref{toy_comp_1}, \fref{toy_comp_2} and \fref{wire_comp} is the faster reduction in the uncertainty for EKLD compared to US.
This observation hints at the faster convergence of the EKLD across all numerical examples.
For expensive problems, with very high-dimensional parameter, space reduced-order model based techniques \cite{bui2008model} need to be used for the context of inferring the statistical expectation of the black-box function. Approaching such problems is beyond the scope of this work.
\begin{figure}[!htbp]
\centering
    \subfigure[]{
        \includegraphics[width=1\columnwidth]{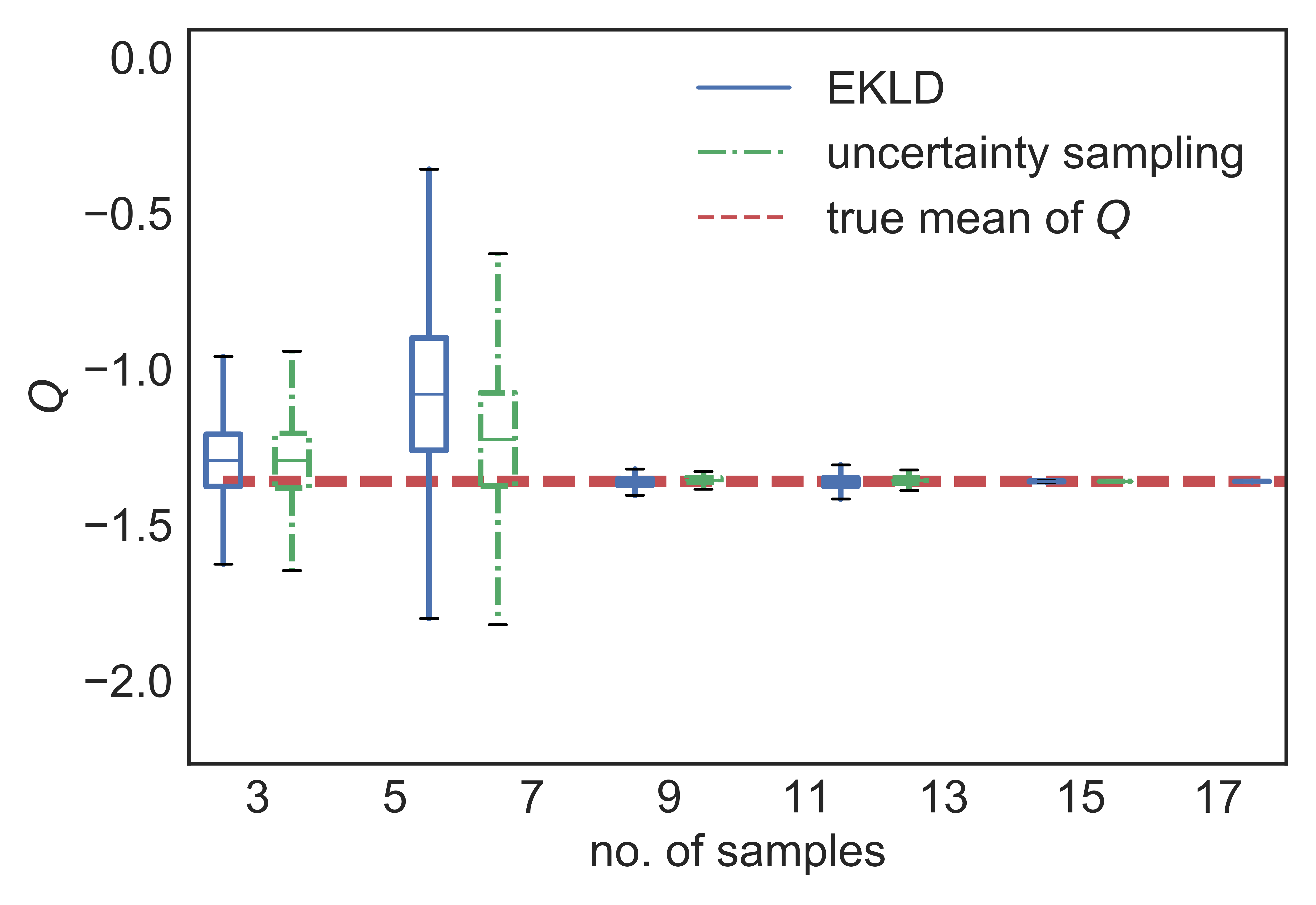}
    }
    \subfigure[]{
        \includegraphics[width=1\columnwidth]{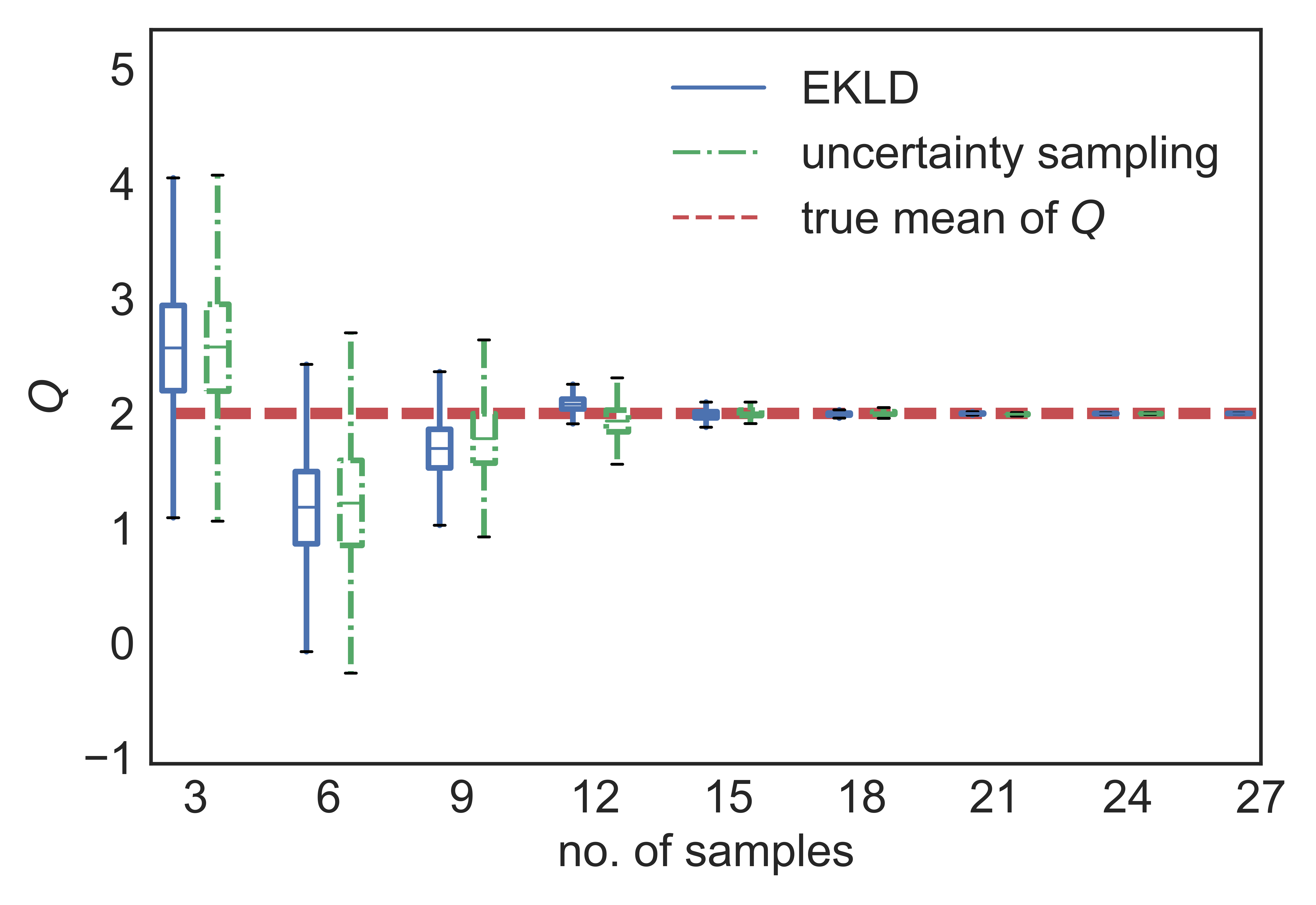}
    }
    \caption{
        Subfigures ~(a), and (b) show the comparison of the EKLD to uncertainty sampling, 
        for synthetic problem nos. 1 and 2 respectively.}
    \label{fig:toy_comp_1}
\end{figure}

\begin{figure}[!htbp]
\centering
    \subfigure[]{
        \includegraphics[width=1\columnwidth]{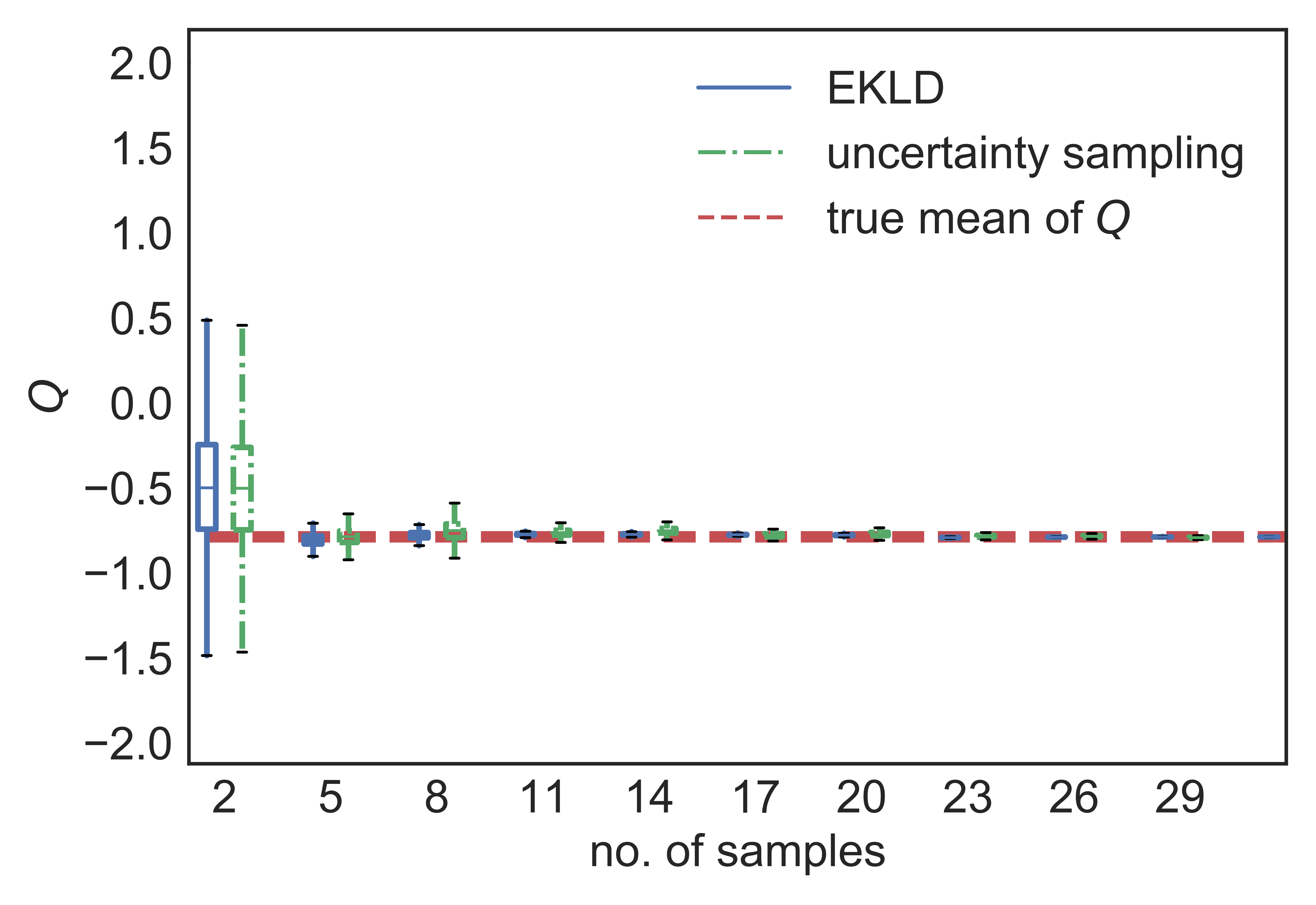}
    }
    \subfigure[]{
        \includegraphics[width=1\columnwidth]{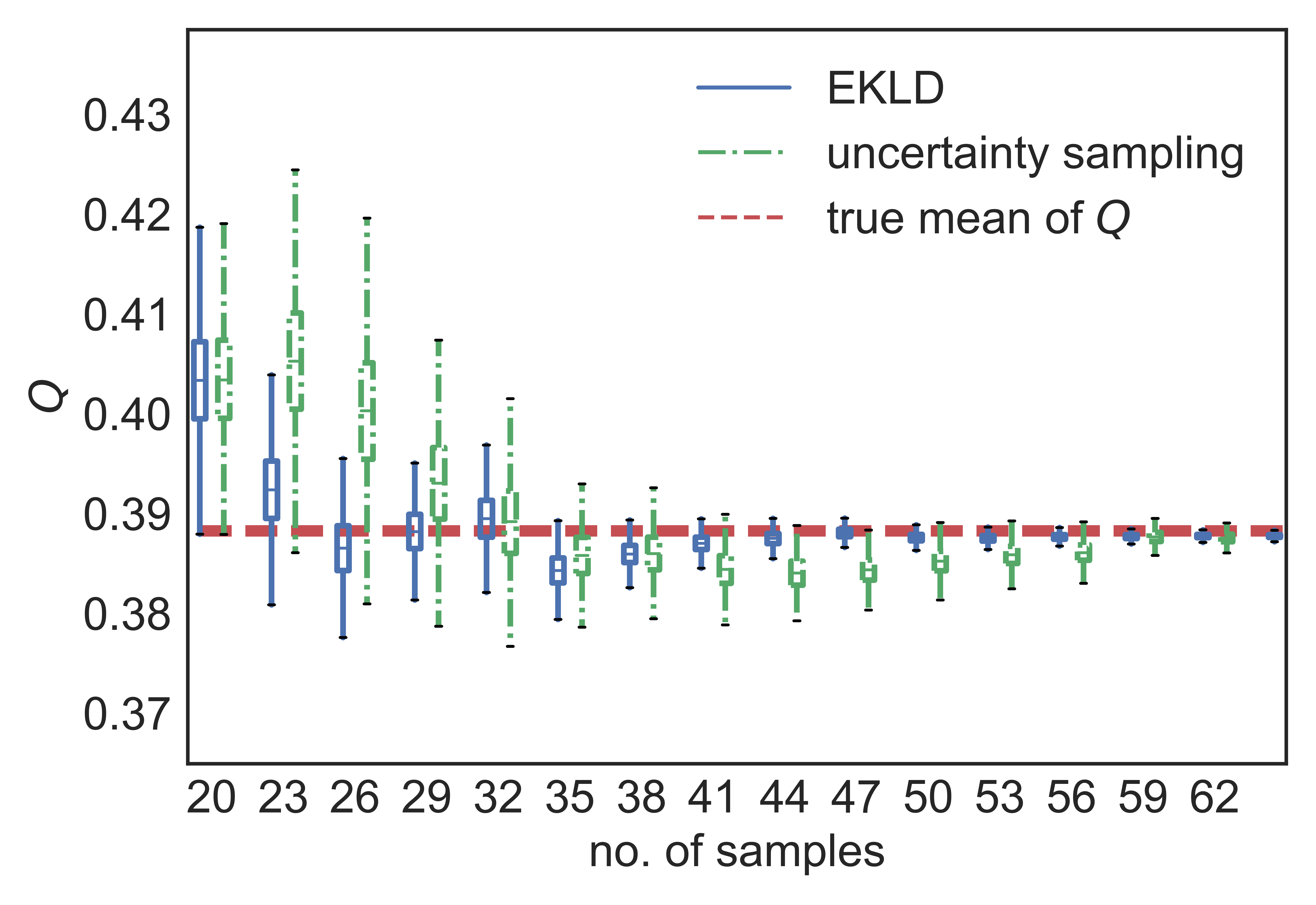}
    }
    \caption{Subfigures ~(a), and (b) show the comparison of the EKLD to uncertainty sampling, 
        for synthetic problem nos. 3 and 4 respectively.}
    \label{fig:toy_comp_2}
\end{figure}
\begin{figure}[!htbp]
\centering
    \includegraphics[width=1\columnwidth]{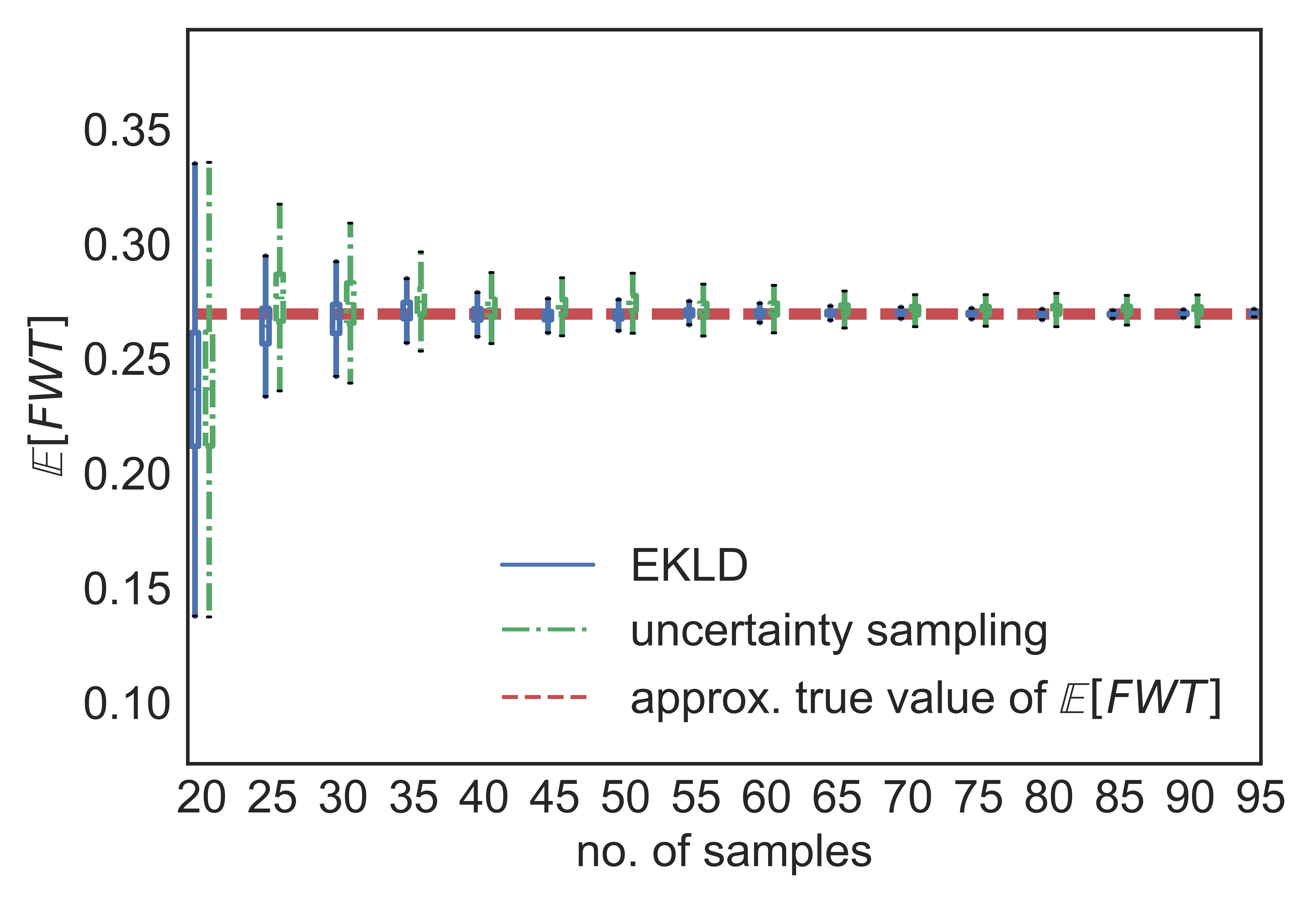}
    \caption{Comparison of the EKLD to uncertainty sampling, 
        for the wire-drawing problem.}
    \label{fig:wire_comp}
\end{figure}
\subsection{Insight into EKLD}
We summarize our thoughts and observations, based on the above experiments, as follows:
\begin{enumerate}
\item{We observe that EKLD and US exhibit similar behavior in low-dimensions, but that EKLD is clearly better in higher-dimensions both in terms of point-wise estimation error and reduction in epistemic uncertainty.
For one-dimensional problems, US took fewer samples to converge to the true value in one of the numerical examples.}
\item{The EKLD quantifies the information gain in the statistical expectation whereas US quantifies the information gain in the parameters (design variables in this work) while selecting the most informative experiment. More work needs to be done to truly analyze and point out the difference between the two methodologies. 
The use of non-stationary GPs is a natural way to fully test the merits and pitfalls of the two methodologies, as it would results in locally adapted designs.}
\item{The number of initial data points differs for each of the above toy problems. This is done on purpose to test the limits of the methodology for examples of varying dimensionality. Thus far the experiments do not reveal a concrete rule for choosing the number of initial data points. However, starting the methodology with too few points can lead to delayed convergence. As a rule of thumb \emph{5d} number of initial points would be considered enough to start the methodology.}
\item{It is also observed that the MCMC samples needed to approximate the EKLD using sample averaging can cause numerical issues. If the MCMC samples are selected from a very short ensemble of chains, the EKLD will be noisy. This would need a more rigorous treatment, than the AEI-based BGO, to optimize the EKLD. To circumvent this issue we do not start the MCMC from scratch at iteration. Instead we use the last particle of the trace from the previous iteration to initialize the MCMC for a given iteration. This results in shorter thermalization times for the MCMC.}
\item{The MCMC details for each problem are in similar vein. The results presented above mention the number of chains and the number of steps per chain for each problem. We observe that the \emph{emcee} \cite{foreman2013emcee} MCMC sampler performs well consistently with a reasonable number of chains and number of steps per chain. One of the requirements of the \emph{emcee} sampler is that the number of chains should be greater than or equal to twice the number of hyper-parameters of the GP model. Thus, the number of chains grows as the dimensionality of the problem increases leading to increased computational cost.}
\end{enumerate}

\section{Conclusions}
\label{sec:conc}
We presented a methodology for designing experiments to infer the value of a particular QoI, the statistical expectation of a physical response.
The methodology leverages the expected KL divergence to compute the information gain in the QoI, from a hypothetical design. 
This work is different from previous work done in sequential design of experiments using KL divergence as it quantifies the information gain in the QoI, instead of the information gain in the model parameters. 
The analytical tractability of the final expressions derived for the expected KL divergence, for learning the statistical expectation of a physical response, obviates computational hurdles induced by sample averaging.

One weakness of our methodology is the assumption that the covariance function of the GP model is stationary. 
The modeling of the hyperparameters of the GP should instead be based on a non-stationary covariance function for more locally adapted designs.
However, the problem of implementing a non-stationary GP is not trivial. 
Another area of limited research is the selection of  number of initial data points, i.e., before starting sequential design of experiments. 
A vast majority of literature on BODE uses \emph{ad hoc} criteria for selection of this initial DOE. 
We accept that this is an open problem and more work is needed in this direction to ensure optimal allocation of budget.
In similar vein, the methodology can be well extended to design experiments to infer generic statistics or quantities of interest which depend on a noisy black-box function. 
We plan to address these challenges in our future work.

\section{Acknowledgments}
This work has been made possible by the financial support provided by National Science Foundation through 
Grant 1662230. 
The authors thank their collaborators at TRDDC, Tata Consultancy Services, Pune, India, 
for providing the steel wire manufacturing problem.
\nolinenumbers
\bibliography{references}
\bibliographystyle{abbrv}
\end{document}